\documentclass{elsarticle}

\usepackage{newtxtext}
\usepackage{textcomp}
\usepackage{stmaryrd}
\usepackage{amsmath}
\usepackage{color} 
\usepackage{amssymb}
\usepackage{amsthm}
\usepackage{MnSymbol}
\usepackage{pgf,tikz}
\usepackage{caption}
\usepackage{accents}
\usepackage{dutchcal} 
\usepackage{chngcntr}
\usepackage{lipsum}
\usepackage{color}
\usepackage[utf8]{inputenc}
\usepackage{pstricks-add}
\usepackage{pgfplots}
\usepackage{centernot}
\usepackage{mathtools}
\usepackage{stmaryrd}
\usepackage{skak}
\usepackage{shorttoc}
\usepackage{wrapfig}
\usepackage{float}
\graphicspath{{./Figures/}}
\usepackage{graphicx}
\usepackage{setspace}
\pgfplotsset{compat=1.15}
\usepackage{mathrsfs}
\usetikzlibrary{arrows}

\usepackage{hyperref, cleveref}
\usepackage[shortlabels]{enumitem}

\numberwithin{equation}{section}

\usetikzlibrary{arrows}
\theoremstyle{plain}
\newtheorem{theorem}{Theorem}[section]

\newtheorem{cor}{Corollary}[section]

\newtheorem{Fact}{Fact}[section]

\newtheorem{notation}{Notation}[section]
\newtheorem{Claim}{Claim}

\numberwithin{equation}{section}
\newtheorem*{convention}{Convention}

\makeatletter
\def\ps@pprintTitle{%
	\let\@oddhead\@empty
	\let\@evenhead\@empty
	\def\@oddfoot{}%
	\let\@evenfoot\@oddfoot}
\makeatother

\begin{document}

	\begin{frontmatter}
		
		\title{Irreducible pairings and indecomposable tournaments}

		\author[add1,add2]{Houmem Belkhechine}
		\ead{houmem.belkhechine@ipeib.ucar.tn}
		\author[add2,add3]{Cherifa Ben Salha}
		\ead{cherifa.bensalha@fsb.u-carthage.tn}
		\author[add2]{Rim Romdhane}
		\ead{rim.romdhane@fsb.u-carthage.tn}
		
		\address[add1]{University of Carthage, Bizerte Preparatory Engineering Institute, 7021, Bizerte, Tunisia}
		\address[add2]{University of Carthage, Faculty of Sciences of Bizerte, UR17ES21 Syst\`{e}mes dynamiques et leurs applications, 7021, Bizerte, Tunisia}
		\address[add3]{University of Tunis, Higher Institute of Applied Studies in Humanities of Zaghouan, 1121, Zaghouan, Tunisia}
		
		\begin{abstract}
			We only consider finite structures. With every totally ordered set $V$ and a subset $P$ of  $\binom{V}{2}$, we associate the underlying tournament ${\rm Inv}(\underline{V}, P)$ obtained from the transitive tournament $\underline{V}:=(V, \{(x,y) \in V \times V : x < y \})$ by reversing $P$, i.e., by reversing the arcs $(x,y)$ such that $\{x,y\} \in P$. The subset $P$ is a pairing (of $\cup P$) if $|\cup P| = 2|P|$, a quasi-pairing (of $\cup P$) if $|\cup P| = 2|P|-1$; it is irreducible if no nontrivial interval of $\cup P$ is a union of connected components of the graph $(\cup P, P)$. In this paper, we consider pairings and quasi-pairings in relation to tournaments. We establish close relationships between irreducibility of pairings (or quasi-pairings) and indecomposability of their underlying tournaments under modular decomposition. For example, given a pairing $P$ of a totally ordered set $V$ of size at least $6$, the pairing $P$ is irreducible if and only if the tournament ${\rm Inv}(\underline{V}, P)$ is indecomposable. This is a consequence of a more general result characterizing indecomposable tournaments obtained from transitive tournaments by reversing pairings. We obtain analogous results in the case of quasi-pairings.   
		\end{abstract}
		\begin{keyword}
			Module \sep pairing  \sep quasi-pairing  \sep indecomposable \sep irreducible \sep transversal 
			\MSC[2020] 05A18 \sep 05C20 \sep 05C35 \sep 06A05 \sep 05D15. 
		\end{keyword}
	\end{frontmatter}
	
	\section{Introduction}
 The structures we consider are all finite. They mainly consist of tournaments, total orders, pairings, and quasi-pairings. We identify a totally ordered set $V$ with the transitive tournament $(V, \{(x,y) \in V \times V : x < y\})$, and we simply denote it by $\underline{V}$ when the total order relation $\leq$ on $V$ is implicitly understood. In this context, $V$ and $\underline{V}$ are often used interchangeably. When $V$ has even size, a pairing of $V$ is a partition $P$ of $V$ whose blocks have size $2$. The pairing $P$ is irreducible if no nontrivial interval of $V$ is a union of blocks. The notion of irreducible pairing was studied by various authors under different considerations (see e.g., \cite{J.S, W.N.H, kleitman, Lehner, A.Nijenhuis, Rstein1, Rstein, Touchard1, Touchard}). However, it seems that irreducible pairings have not been considered in relation to tournaments. In this paper, we address this issue by establishing close relations between irreducibility of pairings (including quasi-parings, see below) and indecomposability of related tournaments under modular decomposition. For more information on indecomposable tournaments, the reader may, for example, refer to \cite{Index, Index2, indecomp1, Cherifa, EFHM, indecomp4, ST}.
 
  With every totally ordered set $V$ and a subset $P$ of $\binom{V}{2}$, we associate the tournament ${\rm Inv}(\underline{V}, P)$ obtained from the transitive tournament $\underline{V}$ by reversing $P$, i.e., by reversing the arcs $(x,y)$ such that $\{x,y\} \in P$. 
 When $P$ is a pairing of a totally ordered set $V$ (of even size at least $6$), $P$ is irreducible if and only if the tournament ${\rm Inv}(\underline{V}, P)$ is indecomposable (see Corollary~\ref{cor V souligné}). This is a consequence of a more general result consisting of a characterization of indecomposable tournaments obtained from $\underline{V}$ by reversing partial pairings, i.e. pairings of subsets of $V$ (see Theorem~\ref{V souligné}).
 
 Since the notion of pairing only involves sets of even sizes, we extend the scope of our study by introducing an analogous notion for sets of odd sizes.
 Given a set $V$ of size $2n+1$, where $n$ is a positive integer, a quasi-pairing of $V$ consists of a family $Q$ of $n+1$ unordered pairs whose union is $V$. When the set $V$ is totally ordered, the quasi-pairing $Q$ is irreducible if no nontrivial interval of $V$ is a union of blocks of the partition of $V$ obtained from $Q$ by merging its two intersecting pairs. Given a quasi-pairing $Q$ of an ordered set $V$ (of odd size at least $7$), let $v_Q^-$ and $v_Q^+$ denote the elements of the symmetric difference of the intersecting pairs of $Q$. The quasi-pairing $Q$ is irreducible if and only if at least one of the tournaments $T$, $T-v_Q^-$, or $T-v_Q^+$ is indecomposable, where $T:={\rm Inv}(\underline{V}, Q)$ (see Corollary~\ref{cor proposition analogue}). This is a consequence of Theorem~\ref{proposition analogue}, which can be seen as an analogue of Theorem~\ref{V souligné} for quasi-pairings. Unlike the case of pairings, the irreducibility of the quasi-pairing $Q$ of $V$ is not a sufficient condition for the tournament ${\rm Inv}(\underline{V}, Q)$ to be indecomposable. Some additional conditions are then required to ensure the indecomposability of ${\rm Inv}(\underline{V}, Q)$ (see Corollary~\ref{cor résultat plus fort}). In a more general setting, these conditions are provided by a second analogue of Theorem~\ref{V souligné} (see Theorem~\ref{résultat plus fort}), which is a characterization of indecomposable tournaments obtained from transitive ones by reversing partial quasi-pairings.      
 
 The original motivation behind our work was to characterize indecomposable tournaments with minimum Slater index. (The Slater index of a tournament is the minimum number of arcs that must be reversed to make it transitive.) This was a question posed by A. Boussa\"{\i}ri during the second author's thesis defense in July 2021. In fact, the indecomposable tournaments with minimum Slater index are obtained from transitive tournaments by reversing some specific irreducible partial pairings or quasi-pairings. They consequently form a subfamily of the family of indecomposable tournaments characterized by Theorems~\ref{V souligné} and \ref{résultat plus fort}. A complete description of these tournaments will soon be provided in a forthcoming paper.

	\section{Basic definitions and preliminaries}
	\subsection{Tournaments and indecomposability}
A {\it  tournament} $T = (V(T), A(T))$ consists of a finite set $V(T)$ of {\it vertices} together with a set $A(T)$ of ordered pairs of distinct vertices, called {\it arcs}, such that for every $x \neq y  \in V(T), (x,y) \in A(T)$ if and only if $(y,x)  \notin A(T)$. 
Given a tournament $T$, the {\it subtournament} of $T$ induced by a subset $X$ of $V(T)$ is the tournament $T[X] := (X, A(T) \cap (X \times X))$. For $X \subseteq V(T)$, the subtournament $T[V(T) \setminus X]$ is also denoted by $T - X$, and by $T - x$ when $X$ is the singleton $\{x\}$. Two tournaments $T$ and $T'$ are {\it isomorphic} if there exists an {\it isomorphism} from $T$ onto $T'$, i.e., a bijection $f$ from $V(T)$ onto $V(T')$ such that for every $x \neq y \in V(T), (x,y) \in A(T)$ if and only if $(f(x), f(y)) \in A(T')$. 
The paper is based on two specific types of tournaments: total orders and indecomposable tournaments.

A {\it total order} is a {\it transitive} tournament, that is, a tournament $T$ such that for every $x, y, z \in V(T)$, if $(x,y) \in A(T)$ and $(y,z) \in A(T)$, then $(x,z) \in A(T)$. We identify a transitive tournament $T$ with the set $V(T)$ totally ordered as follows: for every $x, y \in V(T)$, $x < y$ when $(x,y) \in A(T)$. Given a totally ordered set $V$, when the total ordering $\leq$ on $V$ is implicitly understood, the total order $(V, \{(x,y) \in V \times V : x < y\})$ is denoted by $\underline{V}$. Moreover, $V$ and $\underline{V}$ are often used interchangeably. Since we only consider finite structures, we may assume that $V$ is a subset of $\mathbb{N}$ totally ordered by the natural order on integers. When $V = \{0, \ldots, n-1\}$ for some positive integer $n$, the total order $\underline{V}$ is simply denoted by $\underline{n}$. 

The notion of indecomposability relies on that of a module. The notion of module generalizes to all tournaments the usual notion of interval in a total order. Recall that an {\it interval} of a totally ordered set $V$ is a subset $I$ of $V$ such that every element $v$ in $V \setminus I$ is greater than all the elements of $I$ or smaller than all of them. Analogously, we define a {\it module} of a tournament $T$ to be a subset $M$ of $V(T)$ such that for every vertex $v$ in $V(T) \setminus M$, we have  $\{v\} \times M \subseteq A(T)$ or $M \times \{v\} \subseteq A(T)$. Observe that the notions of module and interval clearly coincide for total orders. The notion of module generalizes also to other combinatorial structures such as graphs and digraphs \cite{indecom}, binary relational structures \cite{ST}, \mbox{2-structures} \cite{ER}, and hypergraphs \cite{Bonizzoni99}. It appears in the literature under various names such as interval \cite{indecom}, convex subset \cite{EFHM}, partitive subset \cite{Sumner}, autonomous set \cite{Buer}, clan \cite{ER}, and, of course, module \cite{CH2}. 

We now come to the notion of indecomposability. Let $T$ be a tournament. The empty set $\emptyset$, the entire vertex set $V(T)$, and its singleton subsets are clearly modules of $T$, called {\it trivial modules}. The tournament $T$ is {\it indecomposable} \cite{Sumner} (or {\it prime} \cite{CH2} or {\it primitive} \cite{ER} or {\it simple} \cite{EFHM}) if all its modules are trivial; otherwise it is {\it decomposable}. Let us consider tournaments with small sizes as examples. 
\pgfplotsset{compat=1.15}
\vspace{-2.02 cm}
\begin{figure}[H]
	\begin{center} 
		\definecolor{wwwwww}{rgb}{0.4,0.4,0.4}
		\begin{tikzpicture}[line cap=round,line join=round,>=triangle 45,x=1cm,y=1cm]
			\clip(-8.30,-0.201538461538482) rectangle (9.62,6.601538461538482);
			\draw [shift={(-0.594262545405895,2.0126639656589345)},line width=0.8pt]  plot[domain=2.436401385820378:3.8299012000075936,variable=\t]({1*1.8460411874569365*cos(\t r)+0*1.8460411874569365*sin(\t r)},{0*1.8460411874569365*cos(\t r)+1*1.8460411874569365*sin(\t r)});
			\draw [shift={(-2.7125095525381355,2.78020276101942)},line width=0.8pt]  plot[domain=-1.025657969332471:1.0421341648568214,variable=\t]({1*1.412648281308104*cos(\t r)+0*1.412648281308104*sin(\t r)},{0*1.412648281308104*cos(\t r)+1*1.412648281308104*sin(\t r)});
			\draw [shift={(-3.405749312678813,2.016618520984345)},line width=0.8pt]  plot[domain=2.5139832212504265:3.7860847908373563,variable=\t]({1*1.9942966417841326*cos(\t r)+0*1.9942966417841326*sin(\t r)},{0*1.9942966417841326*cos(\t r)+1*1.9942966417841326*sin(\t r)});
			\draw [shift={(-3.7680715248076937,2.40923)},line width=0.8pt]  plot[domain=2.2297456528325403:4.053439654347046,variable=\t]({1*2.0120131612118346*cos(\t r)+0*2.0120131612118346*sin(\t r)},{0*2.0120131612118346*cos(\t r)+1*2.0120131612118346*sin(\t r)});
			\draw [shift={(-5.6951087394281075,2.813438380990739)},line width=0.8pt]  plot[domain=-1.0576440992217897:1.0408702989746952,variable=\t]({1*1.37517440178881*cos(\t r)+0*1.37517440178881*sin(\t r)},{0*1.37517440178881*cos(\t r)+1*1.37517440178881*sin(\t r)});
			\draw (0.68,0.6846153846153878) node[anchor=north west] {$W_5$};
			\draw (-2.26,0.683076923076926) node[anchor=north west] {$U_5$};
			\draw (-5.26,0.683076923076926) node[anchor=north west] {$T_5$};
			\draw [shift={(0.8183527896773299,3.2191169638645425)},line width=0.8pt]  plot[domain=-1.3677821399847554:1.3422431904269965,variable=\t]({1*0.8017318910596849*cos(\t r)+0*0.8017318910596849*sin(\t r)},{0*0.8017318910596849*cos(\t r)+1*0.8017318910596849*sin(\t r)});
			\draw [shift={(2.3243741527446793,2.0143403439250487)},line width=0.8pt]  plot[domain=2.424017545002245:3.8760504670855385,variable=\t]({1*1.7844037361103229*cos(\t r)+0*1.7844037361103229*sin(\t r)},{0*1.7844037361103229*cos(\t r)+1*1.7844037361103229*sin(\t r)});
			\draw [line width=0.8pt] (-6,4.383076923076937)-- (-6,0);
			\draw [line width=0.8pt] (-6,0)-- (2,0);
			\draw [line width=0.8pt] (1.98,4.404615384615398)-- (2,0);
			\draw [line width=0.8pt] (1.98,4.404615384615398)-- (-6,4.383076923076937);
			\begin{scriptsize}
				\draw [fill=wwwwww] (-5,4) circle (1.5pt);
				\draw[color=wwwwww] (-4.84,4.1723076923077056) node {0};
				\draw [fill=wwwwww] (-5.02,3.1876923076923167) circle (1.5pt);
				\draw[color=wwwwww] (-4.86,3.3538461538461643) node {1};
				\draw [fill=wwwwww] (-5.02,2.39077) circle (1.5pt);
				\draw[color=wwwwww] (-4.86,2.5784615384615465) node {2};
				\draw [fill=wwwwww] (-5.02,1.6153846153846192) circle (1.5pt);
				\draw[color=wwwwww] (-4.86,1.8815384615384675) node {3};
				\draw [fill=wwwwww] (-5,0.81846) circle (1.5pt);
				\draw[color=wwwwww] (-4.84,0.9846153846153882) node {4};
				\draw [fill=wwwwww] (-2,4) circle (1.5pt);
				\draw[color=wwwwww] (-1.84,4.1723076923077056) node {0};
				\draw [fill=wwwwww] (-2,3.20923) circle (1.5pt);
				\draw[color=wwwwww] (-1.84,3.375384615384626) node {1};
				\draw [fill=wwwwww] (-2,2.455384615384622) circle (1.5pt);
				\draw[color=wwwwww] (-1.84,2.62153846153847) node {2};
				\draw [fill=wwwwww] (-1.98,1.57231) circle (1.5pt);
				\draw[color=wwwwww] (-1.82,1.86) node {3};
				\draw [fill=wwwwww] (-2.02,0.84) circle (1.5pt);
				\draw[color=wwwwww] (-1.86,1.0061538461538497) node {4};
				\draw [fill=wwwwww] (1,4) circle (1.5pt);
				\draw[color=wwwwww] (1.16,4.1723076923077056) node {0};
				\draw [fill=wwwwww] (0.98,3.1876923076923167) circle (1.5pt);
				\draw[color=wwwwww] (1.14,3.3538461538461643) node {1};
				\draw [fill=wwwwww] (0.98,2.43385) circle (1.5pt);
				\draw[color=wwwwww] (1.14,2.65153846153847) node {2};
				\draw [fill=wwwwww] (1,1.6153846153846192) circle (1.5pt);
				\draw[color=wwwwww] (1.16,1.7815384615384675) node {3};
				\draw [fill=wwwwww] (1,0.81846) circle (1.5pt);
				\draw[color=wwwwww] (1.16,0.9846153846153882) node {4};
				\draw [fill=wwwwww,shift={(-2.44,2.04615)}] (0,0) ++(0 pt,2.25pt) -- ++(1.9485571585149868pt,-3.375pt)--++(-3.8971143170299736pt,0 pt) -- ++(1.9485571585149868pt,3.375pt);
				\draw [fill=wwwwww,shift={(-1.3,2.8)}] (0,0) ++(0 pt,2.25pt) -- ++(1.9485571585149868pt,-3.375pt)--++(-3.8971143170299736pt,0 pt) -- ++(1.9485571585149868pt,3.375pt);
				\draw [fill=wwwwww,shift={(-5.4,2.00308)}] (0,0) ++(0 pt,2.25pt) -- ++(1.9485571585149868pt,-3.375pt)--++(-3.8971143170299736pt,0 pt) -- ++(1.9485571585149868pt,3.375pt);
				\draw [fill=wwwwww,shift={(-5.78,2.39077)}] (0,0) ++(0 pt,2.25pt) -- ++(1.9485571585149868pt,-3.375pt)--++(-3.8971143170299736pt,0 pt) -- ++(1.9485571585149868pt,3.375pt);
				\draw [fill=wwwwww,shift={(-4.32,2.8)}] (0,0) ++(0 pt,2.25pt) -- ++(1.9485571585149868pt,-3.375pt)--++(-3.8971143170299736pt,0 pt) -- ++(1.9485571585149868pt,3.375pt);
				\draw [fill=wwwwww,shift={(1.62,3.2307692307692397)}] (0,0) ++(0 pt,2.25pt) -- ++(1.9485571585149868pt,-3.375pt)--++(-3.8971143170299736pt,0 pt) -- ++(1.9485571585149868pt,3.375pt);
				\draw [fill=wwwwww,shift={(0.54,2.02461538461539)}] (0,0) ++(0 pt,2.25pt) -- ++(1.9485571585149868pt,-3.375pt)--++(-3.8971143170299736pt,0 pt) -- ++(1.9485571585149868pt,3.375pt);
			\end{scriptsize}
		\end{tikzpicture}		
		\caption{The indecomposable tournaments on five vertices (missing arcs are oriented from higher to lower).}
		\label{Figure-T-U-W}
	\end{center} 
\end{figure}
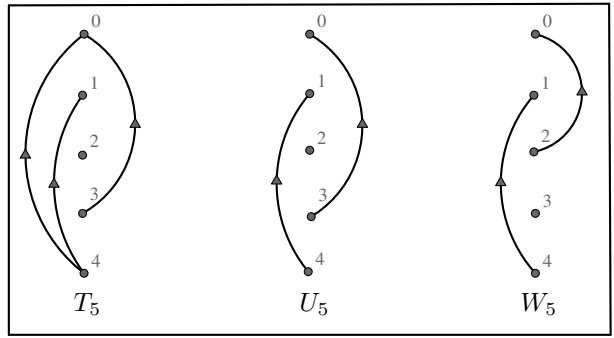 
The tournaments of sizes at most $2$ are clearly indecomposable. The $3$-vertex tournaments are, up to isomorphism, the total order $\underline{3}$, which is decomposable, and the  {\it $3$-cycle} $C_{3} := (\{0,1,2\}, \{(0,1),(1,2),(2,0)\})$, which is indecomposable. Up to isomorphism, there are exactly four $4$-vertex tournaments, all of them are decomposable (see e.g., \cite{Moon}). There are exactly twelve nonisomorphic $5$-vertex tournaments, exactly three of them ($T_5$, $U_5$, and $W_5$; see Figure~\ref{Figure-T-U-W}) are indecomposable (see e.g., \cite{indecomp1}). For sizes at least $5$, it is well-known that there exist indecomposable $n$-vertex tournaments for every $n \geq 5$ (see e.g., \cite{Index}). In fact, Erd\H{o}s et al. \cite{EFHM} proved that almost all tournaments are indecomposable. However, the total orders of sizes at least $3$ are all decomposable.

 \subsection{Reversing arcs, pairings and quasi-pairings} 
 Let $T$ be a tournament and let $P$ be a subset of $\binom{V}{2}$. We denote by ${\rm Inv}(T,P)$ the tournament obtained from $T$ by reversing $P$, i.e., by reversing the arcs $(x,y) \in A(T)$ such that $\{x,y\} \in P$. Thus, ${\rm Inv}(T,P)$ is the tournament defined on $V(T)$ by $A(T) \cap A({\rm Inv}(T,P)) = \{(x,y) \in A(T) : \{x,y\} \notin P\}$. 
 For example, the {\it dual} tournament of $T$ is the tournament $T^\star := {\rm Inv}(T, \binom{V(T)}{2})$. Note that $({\rm Inv}(T,P))^{\star} = {\rm Inv}(T^{\star},P)$. Moreover $T$ and $T^{\star}$ have the same modules; in particular $T$ is indecomposable if and only if $T^{\star}$ is.
These remarks justify that in certain proofs, one may interchange a tournament with its dual.

In this paper, we principally consider tournaments obtained from total orders by reversing pairings or quasi-pairings. The indecomposability of such tournaments is closely related to the irreducibility of the associated pairings or quasi-pairings. Given a totally ordered set $V$, a partition of $V$ (or of $\underline{V}$) is {\it irreducible}  \cite{J.S, E.A} (or {\it connected} \cite{Lehner}) if no nontrivial interval of $V$ is a union of blocks; otherwise it is {\it reducible}. The partitions we need are primarily {\it pairings}, i.e., partitions whose blocks have size $2$.  
An {\it irreducible pairing} \cite{Rstein1, Rstein} is then an irreducible partition whose blocks are unordered pairs. The study of irreducible pairings goes back at least to the 1950s. It seems that Touchard \cite{Touchard1, Touchard} was the first author to consider and study these configurations, which he called proper systems. This name is no longer used, but the notion has been reconsidered by several authors under other names such as irreducible pairings \cite{Rstein1, Rstein}, irreducible diagrams \cite{kleitman}, and linked diagrams \cite{J.S, A.Nijenhuis}.  Given a set $V$ (of even or odd size), a pairing  of a subset (of even size) of $V$ is also called a  {\it partial paring} of $V$. There is a natural one-to-one correspondence between partial pairings of $V$ and involutions without fixed points of subsets (of even sizes) of $V$. To every partial pairing $P$, there corresponds the involution without fixed points of $\cup P$, which we denote by $i_P$, defined by $i_P(B)=B$ for every block $B \in P$. 

We introduce the notion of quasi-pairing as the analogue of that of pairing for sets of odd sizes. Let $W$ be a set of size $2n+1$ with $n \geq 1$. A {\it quasi-pairing} of $W$ is a cover of $W$ by a family of $n+1$ unordered pairs of $W$, i.e., a subset $Q$ of $\binom{W}{2}$ such that $|Q|=n+1$ and $\cup Q =W$. When $W$ is a subset of a set $V$ (of odd or even size), a quasi-pairing of $W$ is also called a {\it partial quasi-pairing} of $V$. We use the following mode of irreducibility for quasi-pairings (see Notation~\ref{quasi-pairing irreducible}).

\begin{notation} \label{not v^v+v-.} \normalfont	
	Let $Q$ be a partial quasi-pairing of a totally ordered set $V$ (of size at least $3$). We denote by $\hat{v}_Q$, $v_Q^-$, and $v_Q^+$ the unique elements of $\cup Q$ such that $\{\hat{v}_Q, v_Q^-\} \in Q$, $\{\hat{v}_Q, v_Q^+\} \in Q$, and $v_Q^- < v_Q^+$. Moreover, the $3$-element set $\{\hat{v}_Q, v_Q^-, v_Q^+\}$ is denoted by $B_Q$. 
\end{notation}

\begin{notation} \label{quasi-pairing irreducible} \normalfont
	Let $Q$ be a quasi-pairing of a totally ordered set $V$ of odd size at least $3$. We denote by $Q_{{\rm part}}$ the partition of $V$ obtained from $Q$ by merging $\{\hat{v}_{Q}, v_{Q}^-\}$ and $\{\hat{v}_{Q}, v_{Q}^+\}$, that is, $Q_{{\rm part}} :=(Q\setminus \{\{\hat{v}_{Q}, v_{Q}^-\}, \{\hat{v}_{Q}, v_{Q}^+\}\})\cup\{B_Q\}$ (see Notation~\ref{not v^v+v-.}). 
	The quasi-pairing $Q$ is {\it irreducible} if the partition $Q_{\rm part}$ is irreducible; otherwise it is {\it reducible}. We also use the following notation.			
	Given $x \in V$, we denote by $\iota_Q(x)$ the vertex subset $\{y \in \cup Q : \{x, y\} \in Q\}$. So $\iota_Q(\hat{v}_{Q}) = \{v_Q^-, v_Q^+\}$, and for every $x \in (\cup Q) \setminus \{\hat{v}_{Q}\}$, we have $|\iota_Q(x)|=1$ and the element of $\iota_Q(x)$ is denoted by $i_Q(x)$.      
\end{notation} 
 
 \subsection{Co-modules and transversals}	
 The notions of co-module and $\Delta$-decomposition were introduced in \cite{Index2} as follows. Given a tournament $T$, a {\it co-module} of $T$ is a subset $M$ of $V(T)$ such that $M$ or $V(T) \setminus M$ is a nontrivial module of $T$. A co-module of $T$ is {\it minimal} if it does not contain any other co-module. We denote by ${\rm mc}(T)$ the family of minimal co-modules of $T$. For example, for every integer $n \geq 3$, we have
 \begin{equation} \label{eq mc(n)}
 	{\rm mc}(\underline{n}) = \{\{0\},\{n-1\}\} \cup \{\{i, i+1 \} : 1 \leq i \leq n-3 \}.
 \end{equation} 
  A {\it co-modular decomposition} of the tournament $T$ is a set of
  pairwise disjoint co-modules of $T$. A {\it $\Delta$-decomposition} of $T$ is a co-modular
  decomposition of $T$ which is of maximum size. Such a size is called the {\it co-modular index} of $T$, and is denoted by $\Delta(T)$. The first two authors \cite{Index2} showed that for every integer $n \geq 3$, we have
  \begin{equation}\label{grand Delta}
  	\Delta(\underline{n}) = \left\lceil \frac{n+1}{2} \right\rceil. 
  \end{equation}
  Note that we need the notions of co-module and $\Delta$-decomposition only in the particular case of total orders, instead of general tournaments.
  
  Given a set family $\mathcal{F}$, a {\it transversal} of $\mathcal{F}$ is any set $R$ that intersects each element of $\mathcal{F}$, that is, such that $F \cap R \neq \varnothing$ for every $F \in \mathcal{F}$.  
  \begin{Fact}\label{inter}
  	Let $V$ be a finite totally ordered set, and let $P$ be a subset of $\binom{V}{2}$. If the tournament ${\rm Inv}(\underline{V}, P)$ is indecomposable, then $\cup P$ is a transversal of the family of all co-modules of $\underline{V}$, and hence a transversal of  ${\rm mc}(\underline{V})$.	
  \end{Fact}
\begin{proof}
	We may assume $V = \llbracket 0, n-1 \rrbracket$ and $\underline{V} = \underline{n}$, where $n$ is a positive integer. Let $M$ be a co-module of $\underline{n}$. If $(\cup P)\cap M = \emptyset$, then $M$ is still a co-module of ${\rm Inv}(\underline{n}, P)$, and hence ${\rm Inv}(\underline{n}, P)$ is decomposable. Therefore, $\cup P$ is a transversal of ${\rm mc}(\underline{V})$. 
	\end{proof}
 
 \section{Presentation of main results}
 Our first result, as presented in Theorem~\ref{V souligné}, provides a characterization of indecomposable tournaments obtained from total orders by reversing (partial) pairings. 
	\begin{theorem} \label{V souligné}
		Let $V$ be a finite totally ordered set such that $|V|\geq 5$, and let $P$ be a partial pairing of $V$. The following assertions are equivalent.
		\begin{enumerate}
			\item The tournament ${\rm Inv}(\underline{V},P)$ is indecomposable.
			\item The partial pairing $P$ of $V$ is an irreducible pairing of a transversal of ${\rm mc}(\underline{V})$.
		\end{enumerate}
	\end{theorem}
For example, the indecomposable tournaments obtained from the total orders $\underline{5}$ and $\underline{6}$ by reversing (partial) pairings are shown in Figure~\ref{Figure-partial-pairing}.
\pgfplotsset{compat=1.15}
\begin{figure}[H]
	\begin{center} 
		\definecolor{wwwwww}{rgb}{0.4,0.4,0.4}
		\begin{tikzpicture}[line cap=round,line join=round,>=triangle 45,x=1cm,y=1cm]
			\clip(-7.94,-3.309230769230792) rectangle (9.98,7.1938461538461755);
			\draw [shift={(-2.4511247009174135,4.812113627685816)},line width=0.8pt]  plot[domain=2.4873410758353223:3.7788177870784203,variable=\t]({1*1.9519448572226428*cos(\t r)+0*1.9519448572226428*sin(\t r)},{0*1.9519448572226428*cos(\t r)+1*1.9519448572226428*sin(\t r)});
			\draw [shift={(-0.562936184955389,4.823721524622741)},line width=0.8pt]  plot[domain=2.45565518442166:3.810346126737798,variable=\t]({1*1.8570900517062767*cos(\t r)+0*1.8570900517062767*sin(\t r)},{0*1.8570900517062767*cos(\t r)+1*1.8570900517062767*sin(\t r)});
			\draw [shift={(1.5234085266153694,4.82538461538462)},line width=0.8pt]  plot[domain=2.484758701219496:3.79842660596009,variable=\t]({1*1.9236670295920621*cos(\t r)+0*1.9236670295920621*sin(\t r)},{0*1.9236670295920621*cos(\t r)+1*1.9236670295920621*sin(\t r)});
			\draw [shift={(3.498727422493864,4.801887258190429)},line width=0.8pt]  plot[domain=2.4672054432228037:3.7991081087288396,variable=\t]({1*1.9187647143465898*cos(\t r)+0*1.9187647143465898*sin(\t r)},{0*1.9187647143465898*cos(\t r)+1*1.9187647143465898*sin(\t r)});
			\draw [shift={(-2.954002523083912,-0.3885000530063606)},line width=0.8pt]  plot[domain=2.667532410049538:3.628213118476708,variable=\t]({1*3.4460169338730298*cos(\t r)+0*3.4460169338730298*sin(\t r)},{0*3.4460169338730298*cos(\t r)+1*3.4460169338730298*sin(\t r)});
			\draw [shift={(-0.25233222376151526,-1.6187689044814755)},line width=0.8pt]  plot[domain=2.518359857222653:3.7485351191283645,variable=\t]({1*2.1276805790848097*cos(\t r)+0*2.1276805790848097*sin(\t r)},{0*2.1276805790848097*cos(\t r)+1*2.1276805790848097*sin(\t r)});
			\draw [shift={(1.4153254437869929,-0.005384615384616616)},line width=0.8pt]  plot[domain=2.451448791275849:3.8317365159037373,variable=\t]({1*1.8353333426641503*cos(\t r)+0*1.8353333426641503*sin(\t r)},{0*1.8353333426641503*cos(\t r)+1*1.8353333426641503*sin(\t r)});
			\draw [shift={(3.527187153846158,-0.04846076923076789)},line width=0.8pt]  plot[domain=2.4921522947526764:3.79103301242691,variable=\t]({1*1.9677839679222373*cos(\t r)+0*1.9677839679222373*sin(\t r)},{0*1.9677839679222373*cos(\t r)+1*1.9677839679222373*sin(\t r)});
			\draw [shift={(-6.822365723961992,4.034370561139446)},line width=0.8pt]  plot[domain=-0.9903954531398806:0.973960956483549,variable=\t]({1*1.4632119081697694*cos(\t r)+0*1.4632119081697694*sin(\t r)},{0*1.4632119081697694*cos(\t r)+1*1.4632119081697694*sin(\t r)});
			\draw [shift={(-2.6075471009976945,4.011787351086579)},line width=0.8pt]  plot[domain=-1.0875998556769249:1.104795227597233,variable=\t]({1*1.3076438310840246*cos(\t r)+0*1.3076438310840246*sin(\t r)},{0*1.3076438310840246*cos(\t r)+1*1.3076438310840246*sin(\t r)});
			\draw [shift={(-1.6095911780527334,3.600769230769238)},line width=0.8pt]  plot[domain=-0.7826502086651228:0.7826502086651231,variable=\t]({1*2.270076186087755*cos(\t r)+0*2.270076186087755*sin(\t r)},{0*2.270076186087755*cos(\t r)+1*2.270076186087755*sin(\t r)});
			\draw [shift={(1.2168565447373139,3.2130749999999995)},line width=0.8pt]  plot[domain=-0.997530809900093:0.9975308099000934,variable=\t]({1*1.4439060312727348*cos(\t r)+0*1.4439060312727348*sin(\t r)},{0*1.4439060312727348*cos(\t r)+1*1.4439060312727348*sin(\t r)});
			\draw [shift={(-6.796001178803625,-1.5723050000000003)},line width=0.8pt]  plot[domain=-0.9833279795734429:0.9833279795734432,variable=\t]({1*1.4361626561367486*cos(\t r)+0*1.4361626561367486*sin(\t r)},{0*1.4361626561367486*cos(\t r)+1*1.4361626561367486*sin(\t r)});
			\draw [shift={(-3.2344725018520832,0.3941709663674229)},line width=0.8pt]  plot[domain=2.3269595007427117:3.9306113929716218,variable=\t]({1*1.0865604579746377*cos(\t r)+0*1.0865604579746377*sin(\t r)},{0*1.0865604579746377*cos(\t r)+1*1.0865604579746377*sin(\t r)});
			\draw [shift={(-3.195943716719675,-2.0030750000000004)},line width=0.8pt]  plot[domain=2.367451187558477:3.9157341196211086,variable=\t]({1*1.1245204269854643*cos(\t r)+0*1.1245204269854643*sin(\t r)},{0*1.1245204269854643*cos(\t r)+1*1.1245204269854643*sin(\t r)});
			\draw [shift={(-5.284094126348764,5.223845)},line width=0.8pt]  plot[domain=2.3158366110762256:3.9673486961033606,variable=\t]({1*1.0559061530047733*cos(\t r)+0*1.0559061530047733*sin(\t r)},{0*1.0559061530047733*cos(\t r)+1*1.0559061530047733*sin(\t r)});
			\draw [shift={(-4.34477007675443,3.6849096622237782)},line width=0.8pt]  plot[domain=-1.2107286597868292:1.1354745830419288,variable=\t]({1*0.8649945671186835*cos(\t r)+0*0.8649945671186835*sin(\t r)},{0*0.8649945671186835*cos(\t r)+1*0.8649945671186835*sin(\t r)});
			\draw [shift={(-3.2906849388461565,-0.41846153846153716)},line width=0.8pt]  plot[domain=-0.8863183823599616:0.8863183823599617,variable=\t]({1*2.041355313680045*cos(\t r)+0*2.041355313680045*sin(\t r)},{0*2.041355313680045*cos(\t r)+1*2.041355313680045*sin(\t r)});
			\draw [shift={(-4.58229548912626,-0.77923)},line width=0.8pt]  plot[domain=-1.1257253101821227:1.1257253101821225,variable=\t]({1*1.352533707364364*cos(\t r)+0*1.352533707364364*sin(\t r)},{0*1.352533707364364*cos(\t r)+1*1.352533707364364*sin(\t r)});
			\draw [shift={(-2.377003195482949,-0.35999024133993723)},line width=0.8pt]  plot[domain=-1.126677501885995:1.1517726493318323,variable=\t]({1*0.8774406016178502*cos(\t r)+0*0.8774406016178502*sin(\t r)},{0*0.8774406016178502*cos(\t r)+1*0.8774406016178502*sin(\t r)});
			\draw [shift={(-0.3553117888029113,-1.1992307692307704)},line width=0.8pt]  plot[domain=-1.1531828826123682:1.1531828826123682,variable=\t]({1*0.8760581191959058*cos(\t r)+0*0.8760581191959058*sin(\t r)},{0*0.8760581191959058*cos(\t r)+1*0.8760581191959058*sin(\t r)});
			\draw [shift={(-1.1599942271482275,-1.1846153846153897)},line width=0.8pt]  plot[domain=-0.9448587222069955:0.9448587222069953,variable=\t]({1*1.9799942271536068*cos(\t r)+0*1.9799942271536068*sin(\t r)},{0*1.9799942271536068*cos(\t r)+1*1.9799942271536068*sin(\t r)});
			\draw [shift={(3.5480135207142864,-1.6369249999999997)},line width=0.8pt]  plot[domain=2.475372369930111:3.807812937249475,variable=\t]({1*1.9690739742171293*cos(\t r)+0*1.9690739742171293*sin(\t r)},{0*1.9690739742171293*cos(\t r)+1*1.9690739742171293*sin(\t r)});
			\draw [shift={(1.2128256296165643,-0.7855970369150697)},line width=0.8pt]  plot[domain=-0.9956884702937732:1.0120712166556172,variable=\t]({1*1.4472104360244291*cos(\t r)+0*1.4472104360244291*sin(\t r)},{0*1.4472104360244291*cos(\t r)+1*1.4472104360244291*sin(\t r)});
			\draw [line width=0.8pt] (-7,6.806153846153866)-- (-7.02,-3.187692307692319);
			\draw [line width=0.8pt] (-7.02,-3.187692307692319)-- (2.98,-3.2092307692307807);
			\draw [line width=0.8pt] (3,6.806153846153866)-- (2.98,-3.2092307692307807);
			\draw [line width=0.8pt] (-7,6.806153846153866)-- (3,6.806153846153866);
			\draw [line width=0.8pt] (-5,6.806153846153866)-- (-5.019916497428793,-3.192000179851703);
			\draw [line width=0.8pt] (-3.04,6.806153846153866)-- (-3.020064946444272,-3.1963075524230544);
			\draw [line width=0.8pt] (-1,6.806153846153866)-- (-0.9999815366537002,-3.2006585013056803);
			\draw [line width=1.2pt] (1.02,6.806153846153866)-- (0.9799164974287926,-3.2049228970713965);
			\draw [line width=0.8pt] (-7.010344786156088,1.6369437792352775)-- (2.9895914393558947,1.593866938990301);
			\begin{scriptsize}
				\draw [fill=wwwwww] (-6,6) circle (1.5pt);
				\draw[color=wwwwww] (-5.84,6.186153846153865) node {0};
				\draw [fill=wwwwww] (-6,5.24462) circle (1.5pt);
				\draw[color=wwwwww] (-5.84,5.432307692307709) node {1};
				\draw [fill=wwwwww] (-6,4.44769) circle (1.5pt);
				\draw[color=wwwwww] (-5.84,4.6138461538461685) node {2};
				\draw [fill=wwwwww] (-6.02,3.60769) circle (1.5pt);
				\draw[color=wwwwww] (-5.86,3.773846153846165) node {3};
				\draw [fill=wwwwww] (-6.02,2.81077) circle (1.5pt);
				\draw[color=wwwwww] (-5.86,3.0984615384615473) node {4};
				\draw [fill=wwwwww] (-4,6) circle (1.5pt);
				\draw[color=wwwwww] (-3.84,6.186153846153865) node {0};
				\draw [fill=wwwwww] (-4,5.24462) circle (1.5pt);
				\draw[color=wwwwww] (-3.84,5.432307692307709) node {1};
				\draw [fill=wwwwww] (-3.98,4.46923) circle (1.5pt);
				\draw[color=wwwwww] (-3.82,4.63538461538463) node {2};
				\draw [fill=wwwwww] (-4.02,3.65077) circle (1.5pt);
				\draw[color=wwwwww] (-3.86,3.83846153846155) node {3};
				\draw [fill=wwwwww] (-4.04,2.875384615384623) circle (1.5pt);
				\draw[color=wwwwww] (-3.88,3.1415384615384704) node {4};
				\draw [fill=wwwwww,shift={(-4.4,4.921538461538476)}] (0,0) ++(0 pt,2.25pt) -- ++(1.9485571585149868pt,-3.375pt)--++(-3.8971143170299736pt,0 pt) -- ++(1.9485571585149868pt,3.375pt);
				\draw [fill=wwwwww] (-2,6) circle (1.5pt);
				\draw[color=wwwwww] (-1.84,6.186153846153865) node {0};
				\draw [fill=wwwwww] (-2.02,5.18) circle (1.5pt);
				\draw[color=wwwwww] (-1.86,5.346153846153863) node {1};
				\draw [fill=wwwwww] (-2.04,4.40462) circle (1.5pt);
				\draw[color=wwwwww] (-1.88,4.592307692307706) node {2};
				\draw [fill=wwwwww] (-2.02,3.67231) circle (1.5pt);
				\draw[color=wwwwww] (-1.86,3.86) node {3};
				\draw [fill=wwwwww] (-2,2.85385) circle (1.5pt);
				\draw[color=wwwwww] (-1.84,3.1415384615384704) node {4};
				\draw [fill=wwwwww] (0,6) circle (1.5pt);
				\draw[color=wwwwww] (0.16,6.186153846153865) node {0};
				\draw [fill=wwwwww] (0,5.201538461538477) circle (1.5pt);
				\draw[color=wwwwww] (0.16,5.367692307692324) node {1};
				\draw [fill=wwwwww] (0,4.4476923076923205) circle (1.5pt);
				\draw[color=wwwwww] (0.16,4.6138461538461685) node {2};
				\draw [fill=wwwwww] (0,3.650769230769241) circle (1.5pt);
				\draw[color=wwwwww] (0.16,3.816923076923089) node {3};
				\draw [fill=wwwwww] (0,2.7892307692307767) circle (1.5pt);
				\draw[color=wwwwww] (0.16,2.955384615384624) node {4};
				\draw [fill=wwwwww] (0,2) circle (1.5pt);
				\draw[color=wwwwww] (0.16,2.38) node {5};
				\draw [fill=wwwwww] (2,6) circle (1.5pt);
				\draw[color=wwwwww] (2.16,6.186153846153865) node {0};
				\draw [fill=wwwwww] (1.98,5.22308) circle (1.5pt);
				\draw[color=wwwwww] (2.14,5.410769230769248) node {1};
				\draw [fill=wwwwww] (2,4.42615) circle (1.5pt);
				\draw[color=wwwwww] (2.16,4.592307692307706) node {2};
				\draw [fill=wwwwww] (1.98,3.62923) circle (1.5pt);
				\draw[color=wwwwww] (2.14,3.795384615384627) node {3};
				\draw [fill=wwwwww] (1.98,2.83231) circle (1.5pt);
				\draw[color=wwwwww] (2.14,3.02) node {4};
				\draw [fill=wwwwww] (2,2) circle (1.5pt);
				\draw[color=wwwwww] (2.16,2.28) node {5};
				\draw [fill=wwwwww] (-6.02,1.184615384615387) circle (1.5pt);
				\draw[color=wwwwww] (-5.86,1.3615384615384653) node {0};
				\draw [fill=wwwwww] (-6,0.46307692307692333) circle (1.5pt);
				\draw[color=wwwwww] (-5.84,0.6292307692307705) node {1};
				\draw [fill=wwwwww] (-6,-0.37692) circle (1.5pt);
				\draw[color=wwwwww] (-5.84,-0.190769230769229587) node {2};
				\draw [fill=wwwwww] (-6,-1.15231) circle (1.5pt);
				\draw[color=wwwwww] (-5.84,-0.9861538461538498) node {3};
				\draw [fill=wwwwww] (-6,-2) circle (1.5pt);
				\draw[color=wwwwww] (-5.84,-1.8261538461538525) node {4};
				\draw [fill=wwwwww] (-6,-2.76769) circle (1.5pt);
				\draw[color=wwwwww] (-5.84,-2.48) node {5};
				\draw [fill=wwwwww] (-3.98,1.184615384615387) circle (1.5pt);
				\draw[color=wwwwww] (-3.82,1.3615384615384653) node {0};
				\draw [fill=wwwwww] (-4,0.44154) circle (1.5pt);
				\draw[color=wwwwww] (-3.84,0.6292307692307705) node {1};
				\draw [fill=wwwwww] (-4,-0.37692) circle (1.5pt);
				\draw[color=wwwwww] (-3.84,-0.190769230769229587) node {2};
				\draw [fill=wwwwww] (-4,-1.21692) circle (1.5pt);
				\draw[color=wwwwww] (-3.84,-1.029230769230773) node {3};
				\draw [fill=wwwwww] (-4,-2) circle (1.5pt);
				\draw[color=wwwwww] (-3.84,-1.7261538461538525) node {4};
				\draw [fill=wwwwww] (-4,-2.78923) circle (1.5pt);
				\draw[color=wwwwww] (-3.84,-2.6015384615384705) node {5};
				\draw [fill=wwwwww] (-2.02,0.44154) circle (1.5pt);
				\draw[color=wwwwww] (-1.86,0.6292307692307705) node {1};
				\draw [fill=wwwwww] (-1.98,-0.37692) circle (1.5pt);
				\draw[color=wwwwww] (-1.82,-0.190769230769229587) node {2};
				\draw [fill=wwwwww] (-2,-1.15231) circle (1.5pt);
				\draw[color=wwwwww] (-1.84,-0.8861538461538498) node {3};
				\draw [fill=wwwwww] (-2,-2) circle (1.5pt);
				\draw[color=wwwwww] (-1.84,-1.6961538461538525) node {4};
				\draw [fill=wwwwww] (-2,-2.83231) circle (1.5pt);
				\draw[color=wwwwww] (-1.84,-2.5861538461538552) node {5};
				\draw [fill=wwwwww] (0,1.1630769230769256) circle (1.5pt);
				\draw[color=wwwwww] (0.16,1.34) node {0};
				\draw [fill=wwwwww] (0,0.42) circle (1.5pt);
				\draw[color=wwwwww] (0.16,0.5861538461538475) node {1};
				\draw [fill=wwwwww] (0,-0.398461538461541) circle (1.5pt);
				\draw[color=wwwwww] (0.16,-0.23230769230769365) node {2};
				\draw [fill=wwwwww] (0,-1.1738461538461589) circle (1.5pt);
				\draw[color=wwwwww] (0.16,-1.0076923076923115) node {3};
				\draw [fill=wwwwww] (0,-2) circle (1.5pt);
				\draw[color=wwwwww] (0.16,-1.7261538461538525) node {4};
				\draw [fill=wwwwww] (0,-2.7892307692307794) circle (1.5pt);
				\draw[color=wwwwww] (0.16,-2.5230769230769322) node {5};
				\draw [fill=wwwwww] (1.96,1.141538461538464) circle (1.5pt);
				\draw[color=wwwwww] (2.12,1.318461538461542) node {0};
				\draw [fill=wwwwww] (1.98,0.4415384615384617) circle (1.5pt);
				\draw[color=wwwwww] (2.14,0.607692307692309) node {1};
				\draw [fill=wwwwww] (2,-0.42) circle (1.5pt);
				\draw[color=wwwwww] (2.16,-0.23230769230769365) node {2};
				\draw [fill=wwwwww] (1.96,-1.23846) circle (1.5pt);
				\draw[color=wwwwww] (2.12,-1.0507692307692347) node {3};
				\draw [fill=wwwwww] (2,-2) circle (1.5pt);
				\draw[color=wwwwww] (2.16,-1.7261538461538525) node {4};
				\draw [fill=wwwwww] (2,-2.85385) circle (1.5pt);
				\draw[color=wwwwww] (2.16,-2.687692307692317) node {5};
				\draw [fill=wwwwww,shift={(-2.42,4.81385)}] (0,0) ++(0 pt,2.25pt) -- ++(1.9485571585149868pt,-3.375pt)--++(-3.8971143170299736pt,0 pt) -- ++(1.9485571585149868pt,3.375pt);
				\draw [fill=wwwwww,shift={(-0.4,4.85692)}] (0,0) ++(0 pt,2.25pt) -- ++(1.9485571585149868pt,-3.375pt)--++(-3.8971143170299736pt,0 pt) -- ++(1.9485571585149868pt,3.375pt);
				\draw [fill=wwwwww,shift={(1.58,4.81385)}] (0,0) ++(0 pt,2.25pt) -- ++(1.9485571585149868pt,-3.375pt)--++(-3.8971143170299736pt,0 pt) -- ++(1.9485571585149868pt,3.375pt);
				\draw [fill=wwwwww,shift={(-6.4,-0.37692)}] (0,0) ++(0 pt,2.25pt) -- ++(1.9485571585149868pt,-3.375pt)--++(-3.8971143170299736pt,0 pt) -- ++(1.9485571585149868pt,3.375pt);
				\draw [fill=wwwwww,shift={(-2.38,-1.62615)}] (0,0) ++(0 pt,2.25pt) -- ++(1.9485571585149868pt,-3.375pt)--++(-3.8971143170299736pt,0 pt) -- ++(1.9485571585149868pt,3.375pt);
				\draw [fill=wwwwww,shift={(-0.42,0)}] (0,0) ++(0 pt,2.25pt) -- ++(1.9485571585149868pt,-3.375pt)--++(-3.8971143170299736pt,0 pt) -- ++(1.9485571585149868pt,3.375pt);
				\draw [fill=wwwwww,shift={(1.56,0)}] (0,0) ++(0 pt,2.25pt) -- ++(1.9485571585149868pt,-3.375pt)--++(-3.8971143170299736pt,0 pt) -- ++(1.9485571585149868pt,3.375pt);
				\draw [fill=wwwwww,shift={(-5.36,3.984615384615396)}] (0,0) ++(0 pt,2.25pt) -- ++(1.9485571585149868pt,-3.375pt)--++(-3.8971143170299736pt,0 pt) -- ++(1.9485571585149868pt,3.375pt);
				\draw [fill=wwwwww,shift={(-1.3,4.02769230769232)}] (0,0) ++(0 pt,2.25pt) -- ++(1.9485571585149868pt,-3.375pt)--++(-3.8971143170299736pt,0 pt) -- ++(1.9485571585149868pt,3.375pt);
				\draw [fill=wwwwww,shift={(0.66,3.553846153846164)}] (0,0) ++(0 pt,2.25pt) -- ++(1.9485571585149868pt,-3.375pt)--++(-3.8971143170299736pt,0 pt) -- ++(1.9485571585149868pt,3.375pt);
				\draw [fill=wwwwww,shift={(2.66,3.166153846153855)}] (0,0) ++(0 pt,2.25pt) -- ++(1.9485571585149868pt,-3.375pt)--++(-3.8971143170299736pt,0 pt) -- ++(1.9485571585149868pt,3.375pt);
				\draw [fill=wwwwww,shift={(-5.36,-1.550769230769237)}] (0,0) ++(0 pt,2.25pt) -- ++(1.9485571585149868pt,-3.375pt)--++(-3.8971143170299736pt,0 pt) -- ++(1.9485571585149868pt,3.375pt);
				\draw [fill=wwwwww,shift={(-4.32,0.4415384615384617)}] (0,0) ++(0 pt,2.25pt) -- ++(1.9485571585149868pt,-3.375pt)--++(-3.8971143170299736pt,0 pt) -- ++(1.9485571585149868pt,3.375pt);
				\draw [fill=wwwwww,shift={(-4.32,-1.9707692307692382)}] (0,0) ++(0 pt,2.25pt) -- ++(1.9485571585149868pt,-3.375pt)--++(-3.8971143170299736pt,0 pt) -- ++(1.9485571585149868pt,3.375pt);
				\draw [fill=wwwwww,shift={(-6.34,5.223076923076938)}] (0,0) ++(0 pt,2.25pt) -- ++(1.9485571585149868pt,-3.375pt)--++(-3.8971143170299736pt,0 pt) -- ++(1.9485571585149868pt,3.375pt);
				\draw [fill=wwwwww,shift={(-3.48,3.704615384615395)}] (0,0) ++(0 pt,2.25pt) -- ++(1.9485571585149868pt,-3.375pt)--++(-3.8971143170299736pt,0 pt) -- ++(1.9485571585149868pt,3.375pt);
				\draw [fill=wwwwww] (-2,1.1630769230769256) circle (1.5pt);
				\draw[color=wwwwww] (-1.84,1.34) node {0};
				\draw [fill=wwwwww,shift={(-1.25,-0.36615)}] (0,0) ++(0 pt,2.25pt) -- ++(1.9485571585149868pt,-3.375pt)--++(-3.8971143170299736pt,0 pt) -- ++(1.9485571585149868pt,3.375pt);
				\draw [fill=wwwwww,shift={(-3.23,-0.7538461538461575)}] (0,0) ++(0 pt,2.25pt) -- ++(1.9485571585149868pt,-3.375pt)--++(-3.8971143170299736pt,0 pt) -- ++(1.9485571585149868pt,3.375pt);
				\draw [fill=wwwwww,shift={(-1.5,-0.38769230769231017)}] (0,0) ++(0 pt,2.25pt) -- ++(1.9485571585149868pt,-3.375pt)--++(-3.8971143170299736pt,0 pt) -- ++(1.9485571585149868pt,3.375pt);
				\draw [fill=wwwwww,shift={(0.52,-1.163076923076928)}] (0,0) ++(0 pt,2.25pt) -- ++(1.9485571585149868pt,-3.375pt)--++(-3.8971143170299736pt,0 pt) -- ++(1.9485571585149868pt,3.375pt);
				\draw [fill=wwwwww,shift={(0.82,-1.18462)}] (0,0) ++(0 pt,2.25pt) -- ++(1.9485571585149868pt,-3.375pt)--++(-3.8971143170299736pt,0 pt) -- ++(1.9485571585149868pt,3.375pt);
				\draw [fill=wwwwww,shift={(1.58,-1.57231)}] (0,0) ++(0 pt,2.25pt) -- ++(1.9485571585149868pt,-3.375pt)--++(-3.8971143170299736pt,0 pt) -- ++(1.9485571585149868pt,3.375pt);
				\draw [fill=wwwwww,shift={(2.66,-0.77538)}] (0,0) ++(0 pt,2.25pt) -- ++(1.9485571585149868pt,-3.375pt)--++(-3.8971143170299736pt,0 pt) -- ++(1.9485571585149868pt,3.375pt);
			\end{scriptsize}
		\end{tikzpicture}
		\caption{Indecomposable tournaments obtained from $\underline{5}$ and $\underline{6}$ by reversing a partial pairing (missing arcs are oriented from higher to lower).}
		\label{Figure-partial-pairing}
	\end{center} 
\end{figure}
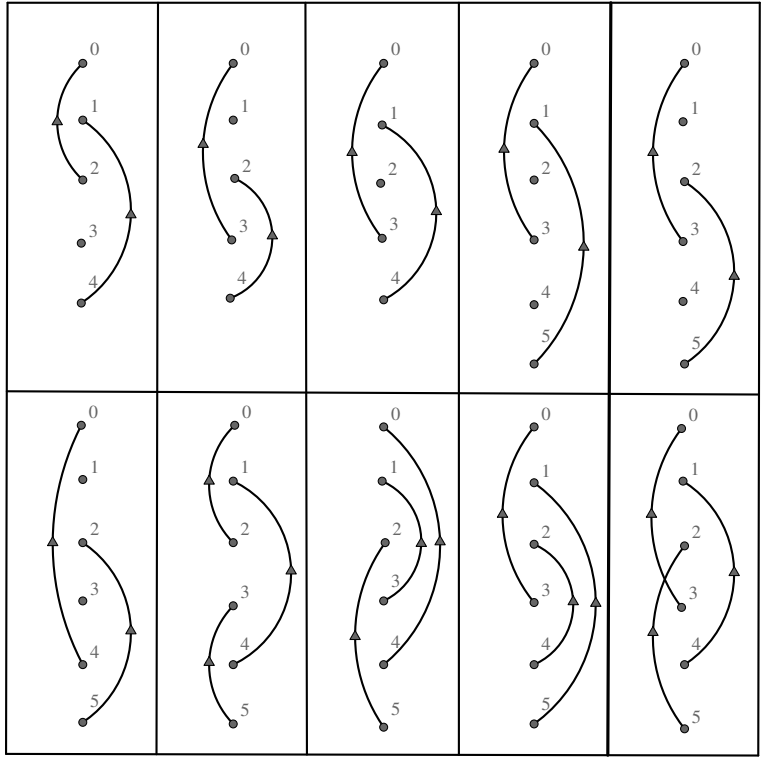
\vspace{-0.5 cm}
The following corollary is a direct consequence of Theorem~\ref{V souligné}. It can be viewed as a novel interpretation of irreducible pairings in terms of indecomposable tournaments.

\begin{cor}\label{cor V souligné} 
	Let $V$ be a finite totally ordered set (with even size) such that $|V|\geq 6$, and let $P$ be a pairing of $V$. The following assertions are equivalent.
	\begin{enumerate}
		\item The pairing $P$ is irreducible.
		\item The tournament ${\rm Inv}(\underline{V},P)$ is indecomposable.
	\end{enumerate}
\end{cor} 
We now provide two analogues of Theorem~\ref{V souligné} for quasi-pairings (see Theorems~\ref{proposition analogue} and \ref{résultat plus fort}).
 
\begin{theorem} \label{proposition analogue}
	Let $V$ be a finite totally ordered set such that $|V|\geq 6$, and let $Q$ be a partial quasi-pairing of $V$. Consider the tournament $T:={\rm Inv}(\underline{V}, Q)$.
	The following two assertions are equivalent.
	\begin{enumerate}
		\item The partial quasi-pairing $Q$ of $V$ is an irreducible quasi-pairing of a transversal of ${\rm mc}(\underline{V})$.
		\item At least one of the tournaments $T$, $T- v_{Q}^-$, or $T- v_{Q}^+$ is indecomposable.
	\end{enumerate}
	Moreover, the second assertion still implies the first one when $|V|=5$.
\end{theorem}
The following corollary, a direct consequence of Theorem~\ref{proposition analogue}, characterizes irreducible quasi-pairings in terms of indecomposable tournaments.
\begin{cor} \label{cor proposition analogue}
	Let $V$ be a finite totally ordered set with odd size such that $|V| \geq 7$. Consider a quasi-pairing $Q$ of $V$, and consider  the tournament $T:={\rm Inv}(\underline{V}, Q)$. The following two assertions are equivalent.
	\begin{enumerate}
		\item The quasi-pairing $Q$ is irreducible.
		\item At least one of the tournaments $T$, $T- v_{Q}^-$, or $T- v_{Q}^+$ is indecomposable.
	\end{enumerate}
\end{cor}
The first condition of Theorem~\ref{proposition analogue} is a necessary but not sufficient condition for the tournament $T$ to be indecomposable. We therefore need some additional conditions to guarantee the indecomposability of $T$. The next theorem provides these conditions optimally. It then characterizes indecomposable tournaments obtained from total orders by reversing (partial) quasi-pairings. 
\begin{theorem} \label{résultat plus fort}
	Let $V$ be a finite totally ordered set such that $|V|=n$ with $n\geq 5$, and let $Q$ be a partial quasi-pairing of $V$. We may assume $V=\llbracket 0, n-1 \rrbracket$ and $\underline{V} = \underline{n}$. Consider the tournament $T:={\rm Inv}(\underline{n}, Q)$.
	The tournament $T$ is indecomposable if and only if the following conditions are satisfied.
	\begin{itemize}
		\item [${\rm (C1)}$] The partial quasi-pairing $Q$ of $\underline{n}$ is an irreducible quasi-pairing of a transversal of ${\rm mc}(\underline{n})$.
		\item [${\rm (C2)}$] $v^+_Q \geq  v^-_Q + 2$.
		\item [${\rm (C3)}$] Given $v \in V(\underline{n})$, if $\{\{v, v+2\}, \{v+1, v+3\}\} \subseteq Q$, then $\hat{v}_Q \in \{v, v+3\}$.
		\item [${\rm (C4)}$] Given $v \in V(\underline{n})$, if $\{v, v+1\} \in Q$, then $\hat{v}_Q \in \{v, v+1\}$ and $\{\hat{v}_Q-1, \hat{v}_Q+1\} \subseteq \cup Q$ (in particular $\hat{v}_Q \notin \{0,n-1\}$).  
	\end{itemize}
\end{theorem}
For example, the indecomposable tournaments obtained from $\underline{5}$ by reversing (partial) quasi-pairings are shown in Figure~\ref{Figure-p-q-p}, where the tournaments $W_5^i$ are isomorphic to $W_5$, and the tournaments $U_5^i$ are isomorphic to $U_5$ (see Figure~\ref{Figure-T-U-W}).

The following corollary is a direct consequence of Theorem~\ref{résultat plus fort}.
\begin{cor}\label{cor résultat plus fort}
	Let $V$ be a finite totally ordered set with odd size such that $|V|=n$ with $n\geq 5$, and let $Q$ be a quasi-pairing of $V$. We may assume $V= \llbracket 0, n-1 \rrbracket$ and $\underline{V}= \underline{n}$. Consider the tournament $T:={\rm Inv}(\underline{n}, Q)$.
The tournament $T$ is indecomposable if and only if the following conditions are satisfied.
\begin{enumerate}
	\item The quasi-pairing $Q$ of $\underline{n}$ is irreducible.
	\item  $v^+_Q \geq  v^-_Q + 2$.
	\item  Given $v \in V(\underline{n})$, if $\{\{v, v+2\}, \{v+1, v+3\}\} \subseteq Q$, then $\hat{v}_Q \in \{v, v+3\}$.
	\item  Given $v \in V(\underline{n})$, if $\{v, v+1\} \in Q$, then $\hat{v}_Q \in \{v, v+1\} \setminus \{0, n-1\}$. 
\end{enumerate}
\end{cor}
The remainder of the paper consists of proofs of Theorems~\ref{V souligné}-\ref{résultat plus fort}.
\pgfplotsset{compat=1.15}
\vspace{-2.5 cm}
\begin{figure}[H] 
	\begin{center} 
	\definecolor{wwwwww}{rgb}{0.4,0.4,0.4}
	\begin{tikzpicture}[line cap=round,line join=round,>=triangle 45,x=1cm,y=1cm]
		\clip(-6.26,-5.266923076923122) rectangle (10.1,7.043076923076947);
		\draw [shift={(-4.480220510000001,3.2169250000000003)},line width=0.8pt]  plot[domain=2.156788911153593:4.126396396025994,variable=\t]({1*0.9398814679786266*cos(\t r)+0*0.9398814679786266*sin(\t r)},{0*0.9398814679786266*cos(\t r)+1*0.9398814679786266*sin(\t r)});
		\draw [shift={(-4.4405693150406575,1.6297796515047296)},line width=0.8pt]  plot[domain=2.1786631875467957:4.079426972186953,variable=\t]({1*0.9795365315307972*cos(\t r)+0*0.9795365315307972*sin(\t r)},{0*0.9795365315307972*cos(\t r)+1*0.9795365315307972*sin(\t r)});
		\draw [shift={(1.0207538651513406,2.39614585603254)},line width=0.8pt]  plot[domain=2.137578534726595:4.120463029511038,variable=\t]({1*1.9011277096352561*cos(\t r)+0*1.9011277096352561*sin(\t r)},{0*1.9011277096352561*cos(\t r)+1*1.9011277096352561*sin(\t r)});
		\draw [shift={(3.3872533022533116,3.2299919524919596)},line width=0.8pt]  plot[domain=2.399957797064482:3.9334100264991165,variable=\t]({1*1.1490280428918884*cos(\t r)+0*1.1490280428918884*sin(\t r)},{0*1.1490280428918884*cos(\t r)+1*1.1490280428918884*sin(\t r)});
		\draw [shift={(6.174810522104649,2.412501199227186)},line width=0.8pt]  plot[domain=2.2078916512140028:4.0626355971331805,variable=\t]({1*1.9749258733691553*cos(\t r)+0*1.9749258733691553*sin(\t r)},{0*1.9749258733691553*cos(\t r)+1*1.9749258733691553*sin(\t r)});
		\draw [shift={(-2.852247157258184,-1.1914259862160586)},line width=0.8pt]  plot[domain=1.9728554637607834:4.360511924041049,variable=\t]({1*0.4286910185766264*cos(\t r)+0*0.4286910185766264*sin(\t r)},{0*0.4286910185766264*cos(\t r)+1*0.4286910185766264*sin(\t r)});
		\draw [shift={(-2.2452356081609097,-2.359629908441109)},line width=0.8pt]  plot[domain=2.3489501607110705:3.960027354542842,variable=\t]({1*1.0752153995178384*cos(\t r)+0*1.0752153995178384*sin(\t r)},{0*1.0752153995178384*cos(\t r)+1*1.0752153995178384*sin(\t r)});
		\draw [shift={(1.2388416937550195,-1.9754592704895775)},line width=0.8pt]  plot[domain=2.1415619590124813:4.190485896352612,variable=\t]({1*0.47909117787979116*cos(\t r)+0*0.47909117787979116*sin(\t r)},{0*0.47909117787979116*cos(\t r)+1*0.47909117787979116*sin(\t r)});
		\draw [shift={(2.0480443584210204,-2.7774385299529922)},line width=0.8pt]  plot[domain=2.295962224071305:4.003699278632631,variable=\t]({1*1.6102960986295747*cos(\t r)+0*1.6102960986295747*sin(\t r)},{0*1.6102960986295747*cos(\t r)+1*1.6102960986295747*sin(\t r)});
		\draw [shift={(3.2001950986008407,-2.00308)},line width=0.8pt]  plot[domain=2.0644354641702236:4.218749843009363,variable=\t]({1*0.46470966672518693*cos(\t r)+0*0.46470966672518693*sin(\t r)},{0*0.46470966672518693*cos(\t r)+1*0.46470966672518693*sin(\t r)});
		\draw [shift={(3.8433451869093913,-3.1954054810522283)},line width=0.8pt]  plot[domain=2.4048972753277145:3.9034805344417975,variable=\t]({1*1.1655914568166807*cos(\t r)+0*1.1655914568166807*sin(\t r)},{0*1.1655914568166807*cos(\t r)+1*1.1655914568166807*sin(\t r)});
		\draw [shift={(-1.9693227962168378,1.6449952608038863)},line width=0.8pt]  plot[domain=2.262140690363164:4.045479491205933,variable=\t]({1*1.0519214883688695*cos(\t r)+0*1.0519214883688695*sin(\t r)},{0*1.0519214883688695*cos(\t r)+1*1.0519214883688695*sin(\t r)});
		\draw [shift={(3.8458619171941266,2.4312626397550248)},line width=0.8pt]  plot[domain=2.263076045297463:4.0327427187285325,variable=\t]({1*2.04586354842406*cos(\t r)+0*2.04586354842406*sin(\t r)},{0*2.04586354842406*cos(\t r)+1*2.04586354842406*sin(\t r)});
		\draw [shift={(-4.87153958174922,-3.590435946356667)},line width=0.8pt]  plot[domain=1.9239623276382338:4.408455317273664,variable=\t]({1*0.42923745537164426*cos(\t r)+0*0.42923745537164426*sin(\t r)},{0*0.42923745537164426*cos(\t r)+1*0.42923745537164426*sin(\t r)});
		\draw [shift={(-4.35223794144737,-2.38)},line width=0.8pt]  plot[domain=2.2616399565530045:4.021545350626582,variable=\t]({1*1.0479834459295838*cos(\t r)+0*1.0479834459295838*sin(\t r)},{0*1.0479834459295838*cos(\t r)+1*1.0479834459295838*sin(\t r)});
		\draw [shift={(-0.8143365863461539,-2.76769)},line width=0.8pt]  plot[domain=2.0261497048893107:4.2570356022902756,variable=\t]({1*0.46765098066373495*cos(\t r)+0*0.46765098066373495*sin(\t r)},{0*0.46765098066373495*cos(\t r)+1*0.46765098066373495*sin(\t r)});
		\draw [shift={(-0.23573433617597078,-1.9987821238038292)},line width=0.8pt]  plot[domain=2.137196386998105:4.129258299105409,variable=\t]({1*1.4242803690055674*cos(\t r)+0*1.4242803690055674*sin(\t r)},{0*1.4242803690055674*cos(\t r)+1*1.4242803690055674*sin(\t r)});
		\draw [shift={(5.203276985043306,-2.7671231313077)},line width=0.8pt]  plot[domain=1.9923849792190964:4.190498016068769,variable=\t]({1*0.44787945714977656*cos(\t r)+0*0.44787945714977656*sin(\t r)},{0*0.44787945714977656*cos(\t r)+1*0.44787945714977656*sin(\t r)});
		\draw [shift={(5.463661112051057,-1.5711154931214546)},line width=0.8pt]  plot[domain=2.135513699521942:4.19924744768028,variable=\t]({1*0.9037403138392142*cos(\t r)+0*0.9037403138392142*sin(\t r)},{0*0.9037403138392142*cos(\t r)+1*0.9037403138392142*sin(\t r)});
		\draw [line width=0.8pt] (-6,4.38308)-- (-6,-5.18);
		\draw [line width=0.8pt] (-6,-5.18)-- (5.98,-5.158461538461557);
		\draw [line width=0.8pt] (5.98,4.404615384615398)-- (5.98,-5.158461538461557);
		\draw [line width=0.8pt] (5.98,4.404615384615398)-- (-6,4.38308);
		\draw [line width=0.8pt] (5.98,-0.39846153846154103)-- (-6,-0.39846153846154103);
		\draw [shift={(4.8056786802940605,2.4338461538461575)},line width=0.8pt]  plot[domain=-1.3609454276919344:1.360945427691934,variable=\t]({1*0.8368197012767176*cos(\t r)+0*0.8368197012767176*sin(\t r)},{0*0.8368197012767176*cos(\t r)+1*0.8368197012767176*sin(\t r)});
		\draw [shift={(-5.221327292419356,-1.6153850000000005)},line width=0.8pt]  plot[domain=-1.3233437524688503:1.3233437524688503,variable=\t]({1*0.8219623679329303*cos(\t r)+0*0.8219623679329303*sin(\t r)},{0*0.8219623679329303*cos(\t r)+1*0.8219623679329303*sin(\t r)});
		\draw [shift={(-3.2389995,-3.1953849999999995)},line width=0.8pt]  plot[domain=-1.28206115013655:1.28206115013655,variable=\t]({1*0.839360506114774*cos(\t r)+0*0.839360506114774*sin(\t r)},{0*0.839360506114774*cos(\t r)+1*0.839360506114774*sin(\t r)});
		\draw [shift={(-1.8333714806940096,-2.7604522958877022)},line width=0.8pt]  plot[domain=-0.9788907668313076:0.9949397901975064,variable=\t]({1*1.4936487323344108*cos(\t r)+0*1.4936487323344108*sin(\t r)},{0*1.4936487323344108*cos(\t r)+1*1.4936487323344108*sin(\t r)});
		\draw [shift={(0.2696459250000001,-1.9923050000000002)},line width=0.8pt]  plot[domain=-1.0346176347441576:1.0346176347441578,variable=\t]({1*1.3905208413015988*cos(\t r)+0*1.3905208413015988*sin(\t r)},{0*1.3905208413015988*cos(\t r)+1*1.3905208413015988*sin(\t r)});
		\draw [shift={(2.2084618827272733,-1.97077)},line width=0.8pt]  plot[domain=-0.9941091557093138:0.9941091557093138,variable=\t]({1*1.4516979291490542*cos(\t r)+0*1.4516979291490542*sin(\t r)},{0*1.4516979291490542*cos(\t r)+1*1.4516979291490542*sin(\t r)});
		\draw [shift={(4.3083462003571436,-2.79154)},line width=0.8pt]  plot[domain=-1.050961034759971:1.0509610347599714,variable=\t]({1*1.3923938200668664*cos(\t r)+0*1.3923938200668664*sin(\t r)},{0*1.3923938200668664*cos(\t r)+1*1.3923938200668664*sin(\t r)});
		\draw (-5.18,0.6184615384615348) node[anchor=north west] {$W_5^1$};
		\draw (-2.96,0.6184615384615348) node[anchor=north west] {$W_5^2$};
		\draw (-0.32,0.6184615384615348) node[anchor=north west] {$W_5^3$};
		\draw (2.32,0.6169230769230732) node[anchor=north west] {$U_5^1$};
		\draw (4.76,0.6169230769230732) node[anchor=north west] {$U_5^2$};
		\draw (-5.22,-4.227692307692333) node[anchor=north west] {$W_5^4$};
		\draw (-3.22,-4.227692307692333) node[anchor=north west] {$W_5^5$};
		\draw (-1.2,-4.227692307692333) node[anchor=north west] {$W_5^6$};
		\draw (0.78,-4.227692307692333) node[anchor=north west] {$W_5^7$};
		\draw (2.8,-4.229230769230794) node[anchor=north west] {$U_5^3$};
		\draw (4.78,-4.227692307692333) node[anchor=north west] {$U_5^4$};
		\draw [shift={(-1.4669160470252636,2.4196050337114947)},line width=0.8pt]  plot[domain=2.207474675276817:4.088258729801959,variable=\t]({1*1.9731353463432515*cos(\t r)+0*1.9731353463432515*sin(\t r)},{0*1.9731353463432515*cos(\t r)+1*1.9731353463432515*sin(\t r)});
		\draw [shift={(0.5665728332753183,3.1913758484162913)},line width=0.8pt]  plot[domain=2.1819669609108563:4.0508414061087965,variable=\t]({1*0.9873589995184598*cos(\t r)+0*0.9873589995184598*sin(\t r)},{0*0.9873589995184598*cos(\t r)+1*0.9873589995184598*sin(\t r)});
		\draw [shift={(2.440403759136786,2.4230769230769265)},line width=0.8pt]  plot[domain=-1.3996536105586017:1.3996536105586022,variable=\t]({1*0.8196669899224147*cos(\t r)+0*0.8196669899224147*sin(\t r)},{0*0.8196669899224147*cos(\t r)+1*0.8196669899224147*sin(\t r)});
		\draw [shift={(5.5993927125506096,1.6261538461538478)},line width=0.8pt]  plot[domain=2.238101498909:4.0450838082705864,variable=\t]({1*1.0008422464021438*cos(\t r)+0*1.0008422464021438*sin(\t r)},{0*1.0008422464021438*cos(\t r)+1*1.0008422464021438*sin(\t r)});
		\draw [line width=0.8pt] (-4.02,-0.39846153846154103)-- (-4.019928953611511,-5.176440093148913);
		\draw [line width=0.8pt] (-2,-0.39846153846154103)-- (-1.9999354828974418,-5.172808411038989);
		\draw [line width=0.8pt] (0.02,-0.39846153846154103)-- (-0.01998060608139518,-5.169248713045487);
		\draw [line width=0.8pt] (1.98,-0.39846153846154103)-- (1.980051652469883,-5.165652918564982);
		\draw [line width=0.8pt] (3.98,-0.39846153846154103)-- (3.980045187830349,-5.1620571937036726);
		\draw [line width=0.8pt] (-3.7999877557735076,4.3870347670983305)-- (-3.82,-0.39846153846154103);
		\draw [line width=0.8pt] (3.820026338645526,4.400732591330223)-- (3.8,-0.39846153846154103);
		\draw [line width=0.8pt] (-1.3999180758067364,4.391349159724528)-- (-1.42,-0.39846153846154103);
		\draw [line width=0.8pt] (1.219996022616984,4.396058747184365)-- (1.24,-0.39846153846154103);
		\begin{scriptsize}
			\draw [fill=wwwwww] (-5,4) circle (1.5pt);
			\draw[color=wwwwww] (-4.84,4.183076923076935) node {0};
			\draw [fill=wwwwww] (-5,3.27385) circle (1.5pt);
			\draw[color=wwwwww] (-4.84,3.4507692307692396) node {1};
			\draw [fill=wwwwww] (-5,2.43385) circle (1.5pt);
			\draw[color=wwwwww] (-4.84,2.610769230769236) node {2};
			\draw [fill=wwwwww] (-5.02,1.65846) circle (1.5pt);
			\draw[color=wwwwww] (-4.86,1.835384615384617) node {3};
			\draw [fill=wwwwww] (-5.02,0.84) circle (1.5pt);
			\draw[color=wwwwww] (-4.86,1.016923076923075) node {4};
			\draw [fill=wwwwww] (-2.64,4.006153846153857) circle (1.5pt);
			\draw[color=wwwwww] (-2.48,4.183076923076935) node {0};
			\draw [fill=wwwwww] (-2.62,3.230769230769238) circle (1.5pt);
			\draw[color=wwwwww] (-2.46,3.4076923076923166) node {1};
			\draw [fill=wwwwww] (-2.64,2.455384615384619) circle (1.5pt);
			\draw[color=wwwwww] (-2.48,2.6323076923076975) node {2};
			\draw [fill=wwwwww] (-2.62,0.8184615384615348) circle (1.5pt);
			\draw[color=wwwwww] (-2.46,0.9953846153846135) node {4};
			\draw [fill=wwwwww] (0,4) circle (1.5pt);
			\draw[color=wwwwww] (0.16,4.183076923076935) node {0};
			\draw [fill=wwwwww] (-0.02,3.209230769230776) circle (1.5pt);
			\draw[color=wwwwww] (0.14,3.386153846153855) node {1};
			\draw [fill=wwwwww] (-0.04,2.412307692307696) circle (1.5pt);
			\draw[color=wwwwww] (0.12,2.5892307692307745) node {2};
			\draw [fill=wwwwww] (-0.02,1.6153846153846152) circle (1.5pt);
			\draw[color=wwwwww] (0.14,1.792307692307694) node {3};
			\draw [fill=wwwwww] (-0.04,0.8184615384615348) circle (1.5pt);
			\draw[color=wwwwww] (0.12,0.9953846153846135) node {4};
			\draw [fill=wwwwww] (-2.64,1.636923076923077) circle (1.5pt);
			\draw[color=wwwwww] (-2.48,1.8138461538461554) node {3};
			\draw [fill=wwwwww] (2.54,4.006153846153857) circle (1.5pt);
			\draw[color=wwwwww] (2.7,4.183076923076935) node {0};
			\draw [fill=wwwwww] (2.58,3.230769230769238) circle (1.5pt);
			\draw[color=wwwwww] (2.74,3.4076923076923166) node {1};
			\draw [fill=wwwwww] (2.58,2.412307692307696) circle (1.5pt);
			\draw[color=wwwwww] (2.74,2.5892307692307745) node {2};
			\draw [fill=wwwwww] (2.58,1.6153846153846152) circle (1.5pt);
			\draw[color=wwwwww] (2.74,1.822307692307694) node {3};
			\draw [fill=wwwwww] (2.56,0.84) circle (1.5pt);
			\draw[color=wwwwww] (2.72,1.016923076923075) node {4};
			\draw [fill=wwwwww] (5,4) circle (1.5pt);
			\draw[color=wwwwww] (5.16,4.183076923076935) node {0};
			\draw [fill=wwwwww] (4.98,3.2523076923076997) circle (1.5pt);
			\draw[color=wwwwww] (5.14,3.429230769230778) node {1};
			\draw [fill=wwwwww] (4.98,2.412307692307696) circle (1.5pt);
			\draw[color=wwwwww] (5.14,2.5892307692307745) node {2};
			\draw [fill=wwwwww] (4.98,1.6153846153846152) circle (1.5pt);
			\draw[color=wwwwww] (5.14,1.822307692307694) node {3};
			\draw [fill=wwwwww] (4.98,0.84) circle (1.5pt);
			\draw[color=wwwwww] (5.14,1.016923076923075) node {4};
			\draw [fill=wwwwww] (-5.02,-0.81846) circle (1.5pt);
			\draw[color=wwwwww] (-4.86,-0.64153846153847087) node {0};
			\draw [fill=wwwwww] (-5.02,-1.57231) circle (1.5pt);
			\draw[color=wwwwww] (-4.86,-1.3953846153846281) node {1};
			\draw [fill=wwwwww] (-5.02,-2.41231) circle (1.5pt);
			\draw[color=wwwwww] (-4.86,-2.135384615384632) node {2};
			\draw [fill=wwwwww] (-5.02,-3.18769) circle (1.5pt);
			\draw[color=wwwwww] (-4.86,-3.010769230769251) node {3};
			\draw [fill=wwwwww] (-5,-4) circle (1.5pt);
			\draw[color=wwwwww] (-4.84,-3.829230769230793) node {4};
			\draw [fill=wwwwww] (-3.02,-0.79692) circle (1.5pt);
			\draw[color=wwwwww] (-2.86,-0.62) node {0};
			\draw [fill=wwwwww] (-3,-1.59385) circle (1.5pt);
			\draw[color=wwwwww] (-2.84,-1.4169230769230897) node {1};
			\draw [fill=wwwwww] (-3,-2.39077) circle (1.5pt);
			\draw[color=wwwwww] (-2.84,-2.21384615384617) node {2};
			\draw [fill=wwwwww] (-2.98,-3.14462) circle (1.5pt);
			\draw[color=wwwwww] (-2.82,-2.9676923076923274) node {3};
			\draw [fill=wwwwww] (-3,-4) circle (1.5pt);
			\draw[color=wwwwww] (-2.84,-3.729230769230793) node {4};
			\draw [fill=wwwwww] (-1,-0.79692) circle (1.5pt);
			\draw[color=wwwwww] (-0.84,-0.62) node {0};
			\draw [fill=wwwwww] (-1.02,-1.50769) circle (1.5pt);
			\draw[color=wwwwww] (-0.86,-1.3307692307692432) node {1};
			\draw [fill=wwwwww] (-1.02,-2.34769) circle (1.5pt);
			\draw[color=wwwwww] (-0.86,-2.170769230769247) node {2};
			\draw [fill=wwwwww] (-1.02,-3.18769) circle (1.5pt);
			\draw[color=wwwwww] (-0.86,-3.010769230769251) node {3};
			\draw [fill=wwwwww] (-1,-4) circle (1.5pt);
			\draw[color=wwwwww] (-0.84,-3.729230769230793) node {4};
			\draw [fill=wwwwww] (0.98,-0.79692) circle (1.5pt);
			\draw[color=wwwwww] (1.14,-0.62) node {0};
			\draw [fill=wwwwww] (0.98,-1.57231) circle (1.5pt);
			\draw[color=wwwwww] (1.14,-1.3953846153846281) node {1};
			\draw [fill=wwwwww] (1,-2.39077) circle (1.5pt);
			\draw[color=wwwwww] (1.16,-2.21384615384617) node {2};
			\draw [fill=wwwwww] (0.98,-3.18769) circle (1.5pt);
			\draw[color=wwwwww] (1.14,-2.910769230769251) node {3};
			\draw [fill=wwwwww] (1,-4) circle (1.5pt);
			\draw[color=wwwwww] (1.16,-3.829230769230793) node {4};
			\draw [fill=wwwwww] (3,-0.75385) circle (1.5pt);
			\draw[color=wwwwww] (3.16,-0.576923076923086) node {0};
			\draw [fill=wwwwww] (2.98,-1.59385) circle (1.5pt);
			\draw[color=wwwwww] (3.14,-1.4169230769230897) node {1};
			\draw [fill=wwwwww] (2.98,-2.41231) circle (1.5pt);
			\draw[color=wwwwww] (3.14,-2.235384615384632) node {2};
			\draw [fill=wwwwww] (3,-3.18769) circle (1.5pt);
			\draw[color=wwwwww] (3.16,-2.910769230769251) node {3};
			\draw [fill=wwwwww] (3,-4) circle (1.5pt);
			\draw[color=wwwwww] (3.16,-3.829230769230793) node {4};
			\draw [fill=wwwwww,shift={(-5.42,3.23077)}] (0,0) ++(0 pt,2.25pt) -- ++(1.9485571585149868pt,-3.375pt)--++(-3.8971143170299736pt,0 pt) -- ++(1.9485571585149868pt,3.375pt);
			\draw [fill=wwwwww,shift={(-5.42,1.61538)}] (0,0) ++(0 pt,2.25pt) -- ++(1.9485571585149868pt,-3.375pt)--++(-3.8971143170299736pt,0 pt) -- ++(1.9485571585149868pt,3.375pt);
			\draw [fill=wwwwww,shift={(-0.88,2.4338461538461575)}] (0,0) ++(0 pt,2.25pt) -- ++(1.9485571585149868pt,-3.375pt)--++(-3.8971143170299736pt,0 pt) -- ++(1.9485571585149868pt,3.375pt);
			\draw [fill=wwwwww,shift={(2.24,3.1661538461538528)}] (0,0) ++(0 pt,2.25pt) -- ++(1.9485571585149868pt,-3.375pt)--++(-3.8971143170299736pt,0 pt) -- ++(1.9485571585149868pt,3.375pt);
			\draw [fill=wwwwww,shift={(4.2,2.4338461538461575)}] (0,0) ++(0 pt,2.25pt) -- ++(1.9485571585149868pt,-3.375pt)--++(-3.8971143170299736pt,0 pt) -- ++(1.9485571585149868pt,3.375pt);
			\draw [fill=wwwwww,shift={(-3.28,-1.16308)}] (0,0) ++(0 pt,2.25pt) -- ++(1.9485571585149868pt,-3.375pt)--++(-3.8971143170299736pt,0 pt) -- ++(1.9485571585149868pt,3.375pt);
			\draw [fill=wwwwww,shift={(-3.32,-2.3907692307692394)}] (0,0) ++(0 pt,2.25pt) -- ++(1.9485571585149868pt,-3.375pt)--++(-3.8971143170299736pt,0 pt) -- ++(1.9485571585149868pt,3.375pt);
			\draw [fill=wwwwww,shift={(0.76,-1.96)}] (0,0) ++(0 pt,2.25pt) -- ++(1.9485571585149868pt,-3.375pt)--++(-3.8971143170299736pt,0 pt) -- ++(1.9485571585149868pt,3.375pt);
			\draw [fill=wwwwww,shift={(0.44,-2.69231)}] (0,0) ++(0 pt,2.25pt) -- ++(1.9485571585149868pt,-3.375pt)--++(-3.8971143170299736pt,0 pt) -- ++(1.9485571585149868pt,3.375pt);
			\draw [fill=wwwwww,shift={(2.74,-1.9384615384615458)}] (0,0) ++(0 pt,2.25pt) -- ++(1.9485571585149868pt,-3.375pt)--++(-3.8971143170299736pt,0 pt) -- ++(1.9485571585149868pt,3.375pt);
			\draw [fill=wwwwww,shift={(2.68,-3.123076923076934)}] (0,0) ++(0 pt,2.25pt) -- ++(1.9485571585149868pt,-3.375pt)--++(-3.8971143170299736pt,0 pt) -- ++(1.9485571585149868pt,3.375pt);
			\draw [fill=wwwwww,shift={(-3.02,1.5938461538461537)}] (0,0) ++(0 pt,2.25pt) -- ++(1.9485571585149868pt,-3.375pt)--++(-3.8971143170299736pt,0 pt) -- ++(1.9485571585149868pt,3.375pt);
			\draw [fill=wwwwww,shift={(1.8,2.4338461538461575)}] (0,0) ++(0 pt,2.25pt) -- ++(1.9485571585149868pt,-3.375pt)--++(-3.8971143170299736pt,0 pt) -- ++(1.9485571585149868pt,3.375pt);
			\draw [fill=wwwwww,shift={(-5.3,-3.56462)}] (0,0) ++(0 pt,2.25pt) -- ++(1.9485571585149868pt,-3.375pt)--++(-3.8971143170299736pt,0 pt) -- ++(1.9485571585149868pt,3.375pt);
			\draw [fill=wwwwww,shift={(-5.4,-2.35846)}] (0,0) ++(0 pt,2.25pt) -- ++(1.9485571585149868pt,-3.375pt)--++(-3.8971143170299736pt,0 pt) -- ++(1.9485571585149868pt,3.375pt);
			\draw [fill=wwwwww,shift={(-1.28,-2.72462)}] (0,0) ++(0 pt,2.25pt) -- ++(1.9485571585149868pt,-3.375pt)--++(-3.8971143170299736pt,0 pt) -- ++(1.9485571585149868pt,3.375pt);
			\draw [fill=wwwwww,shift={(-1.66,-1.99231)}] (0,0) ++(0 pt,2.25pt) -- ++(1.9485571585149868pt,-3.375pt)--++(-3.8971143170299736pt,0 pt) -- ++(1.9485571585149868pt,3.375pt);
			\draw [fill=wwwwww] (4.98,-0.80769) circle (1.5pt);
			\draw[color=wwwwww] (5.14,-0.62) node {0};
			\draw [fill=wwwwww] (5,-1.58308) circle (1.5pt);
			\draw[color=wwwwww] (5.16,-1.4169230769230897) node {1};
			\draw [fill=wwwwww] (5.02,-2.35846) circle (1.5pt);
			\draw[color=wwwwww] (5.18,-2.170769230769247) node {2};
			\draw [fill=wwwwww] (4.98,-3.15538) circle (1.5pt);
			\draw[color=wwwwww] (5.14,-2.9676923076923274) node {3};
			\draw [fill=wwwwww] (5,-4) circle (1.5pt);
			\draw[color=wwwwww] (5.16,-3.729230769230793) node {4};
			\draw [fill=wwwwww,shift={(4.76,-2.70308)}] (0,0) ++(0 pt,2.25pt) -- ++(1.9485571585149868pt,-3.375pt)--++(-3.8971143170299736pt,0 pt) -- ++(1.9485571585149868pt,3.375pt);
			\draw [fill=wwwwww,shift={(4.56,-1.58308)}] (0,0) ++(0 pt,2.25pt) -- ++(1.9485571585149868pt,-3.375pt)--++(-3.8971143170299736pt,0 pt) -- ++(1.9485571585149868pt,3.375pt);
			\draw [fill=wwwwww,shift={(5.64,2.3692307692307724)}] (0,0) ++(0 pt,2.25pt) -- ++(1.9485571585149868pt,-3.375pt)--++(-3.8971143170299736pt,0 pt) -- ++(1.9485571585149868pt,3.375pt);
			\draw [fill=wwwwww,shift={(-4.4,-1.58308)}] (0,0) ++(0 pt,2.25pt) -- ++(1.9485571585149868pt,-3.375pt)--++(-3.8971143170299736pt,0 pt) -- ++(1.9485571585149868pt,3.375pt);
			\draw [fill=wwwwww,shift={(-2.4,-3.22)}] (0,0) ++(0 pt,2.25pt) -- ++(1.9485571585149868pt,-3.375pt)--++(-3.8971143170299736pt,0 pt) -- ++(1.9485571585149868pt,3.375pt);
			\draw [fill=wwwwww,shift={(-0.34,-2.78923)}] (0,0) ++(0 pt,2.25pt) -- ++(1.9485571585149868pt,-3.375pt)--++(-3.8971143170299736pt,0 pt) -- ++(1.9485571585149868pt,3.375pt);
			\draw [fill=wwwwww,shift={(1.66,-1.97077)}] (0,0) ++(0 pt,2.25pt) -- ++(1.9485571585149868pt,-3.375pt)--++(-3.8971143170299736pt,0 pt) -- ++(1.9485571585149868pt,3.375pt);
			\draw [fill=wwwwww,shift={(3.66,-1.99231)}] (0,0) ++(0 pt,2.25pt) -- ++(1.9485571585149868pt,-3.375pt)--++(-3.8971143170299736pt,0 pt) -- ++(1.9485571585149868pt,3.375pt);
			\draw [fill=wwwwww,shift={(5.7,-2.74615)}] (0,0) ++(0 pt,2.25pt) -- ++(1.9485571585149868pt,-3.375pt)--++(-3.8971143170299736pt,0 pt) -- ++(1.9485571585149868pt,3.375pt);
			\draw [fill=wwwwww,shift={(-3.44,2.4338461538461575)}] (0,0) ++(0 pt,2.25pt) -- ++(1.9485571585149868pt,-3.375pt)--++(-3.8971143170299736pt,0 pt) -- ++(1.9485571585149868pt,3.375pt);
			\draw [fill=wwwwww,shift={(-0.42,3.230769230769238)}] (0,0) ++(0 pt,2.25pt) -- ++(1.9485571585149868pt,-3.375pt)--++(-3.8971143170299736pt,0 pt) -- ++(1.9485571585149868pt,3.375pt);
			\draw [fill=wwwwww,shift={(3.26,2.412307692307696)}] (0,0) ++(0 pt,2.25pt) -- ++(1.9485571585149868pt,-3.375pt)--++(-3.8971143170299736pt,0 pt) -- ++(1.9485571585149868pt,3.375pt);
			\draw [fill=wwwwww,shift={(4.6,1.572307692307692)}] (0,0) ++(0 pt,2.25pt) -- ++(1.9485571585149868pt,-3.375pt)--++(-3.8971143170299736pt,0 pt) -- ++(1.9485571585149868pt,3.375pt);
		\end{scriptsize}
	\end{tikzpicture}		
\caption{Indecomposable tournaments obtained from $\underline{5}$ by reversing a partial quasi-pairing (missing arcs are oriented from higher to lower).}
\label{Figure-p-q-p}
\end{center} 
\end{figure}
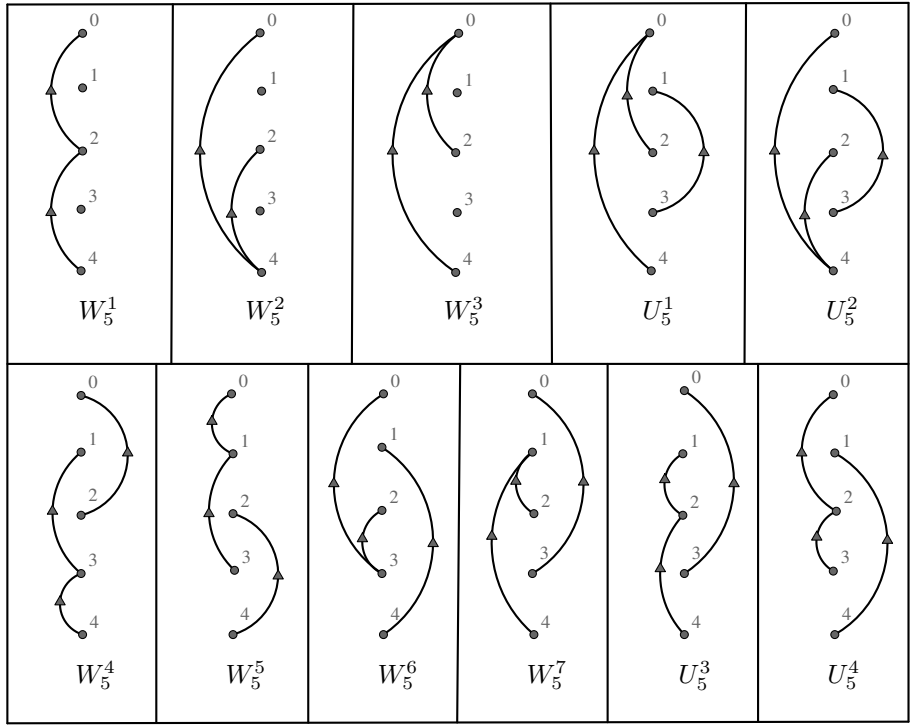 
\vspace{-1 cm} 	
\section{Proofs of Theorems \ref{V souligné}--\ref{résultat plus fort}}\label{sect proof of theorem pairing}
We first introduce some practical notations. Given two nonnegative integers $p$ and $q$, we denote by $\llbracket p, q \rrbracket$, $\llbracket p, q \llbracket$, $\rrbracket p, q \rrbracket$, and $\rrbracket p, q \llbracket$ the intervals of $\mathbb{N}$ defined as follows: $\llbracket p, q \rrbracket := \{i \in \mathbb{N} : p \leq i \leq q\}$; $\llbracket p, q \llbracket := \llbracket p, q \rrbracket \setminus \{q\}$; $\rrbracket p, q \rrbracket := \llbracket p, q \rrbracket \setminus \{p\}$; and $\rrbracket p, q \llbracket := \llbracket p, q \rrbracket \setminus\{p,q\}$.
The next notations are practical when dealing with modules.
\begin{convention} \normalfont
	Let $T$ be a tournament. For $X \subseteq V(T)$, $\overline{X}$ denotes $V(T) \setminus X$.
\end{convention} 
Given a tournament $T$, we also denote by $T$ the function
$$\begin{array}{ccccc}
	T & : & [ V(T) ]^2 & \longrightarrow & \{0, 1\} \\[6pt]
	& & (x, y) & \longmapsto & T(x, y) =  \left \{
	\begin{array}{r c l}
		1  & \text{if} & (x, y) \in A(T), \\[6pt]
		
		0 & \text{if} & (x, y) \notin A(T),
	\end{array}
	\right .\\
\end{array}$$
where $[V(T)]^2 := (V(T)\times V(T))\setminus\{(x, x) : x \in V(T)\}$. 
Moreover, given $X \subseteq V(T)$ and $v \in \overline{X}$, the notation $T(v,X) = i$, where $i \in \{0, 1\}$, signifies that $T(v,x) =i$ for every $x \in X$. The subset $X$ is then a module of $T$ if and only if for every $v \in \overline{X}$, there is $i \in \{0, 1\}$ such that $T(v,X) =i$.
	
	\begin{proof}[Proof of Theorem~\ref{V souligné}]
		We may assume $V = \llbracket 0, n-1 \rrbracket$ and $\underline{V} = \underline{n}$, where $n \geq 5$. Let $T$ and $X$ denote ${\rm Inv}(\underline{n},P)$ and $\cup P$, respectively. Recall that $P$ is a pairing of $X$. To begin, suppose that the tournament $T$ is indecomposable. By Fact~\ref{inter}, $X$ is a transversal of ${\rm mc}(\underline{n})$. Suppose for a contradiction that the pairing $P$ of $\underline{X}$ is reducible. There exists $R \varsubsetneq P$ such that $\cup R$ is a nontrivial interval of $\underline{X}$. Set $M:=\llbracket \min(\cup R), \max(\cup R) \rrbracket$. If there exists $x$ in $M \cap (\cup (P \setminus R))$, then $x \in X \setminus (\cup R)$ and $\min(\cup R) < x < \max(\cup R)$, which contradicts that $\cup R$ is an interval of $\underline{X}$. Therefore $\cup (P \setminus R) \subseteq \overline{M}$. Since $R \subseteq \binom{M}{2}$, it follows that $P \subseteq \binom{\overline{M}}{2} \cup \binom{M}{2}$ and hence $M$ is a module of $T$. Moreover, $|M| \geq 2$ because $|\cup R| \geq 2$, and $|\overline{M}| \geq 2$ because $\cup (P \setminus R) \subseteq \overline{M}$ and $|\cup (P \setminus R)| \geq 2$. Therefore, the module $M$ of $T$ is nontrivial, contradicting that $T$ is indecomposable. So the second assertion holds. 
		
		Conversely, suppose that the tournament $T$ is decomposable. Suppose that
		\begin{equation} \label{eqco} 
			X \text{ is a transversal of } {\rm mc}(\underline{n}) \ (\text{see } (\ref{eq mc(n)})).
		\end{equation}
		We have to prove that the pairing $P$ of $\underline{X}$ is reducible.  
		Since $\Delta(\underline{n}) = \left\lceil \frac{n+1}{2} \right\rceil \geq 3$ (see (\ref{grand Delta})), and since $|X|$ is even because $P$ is a pairing, it follows from (\ref{eqco}) that $|X| \geq 4$. Now consider a nontrivial module $M$ of $T$, and let $m^-$ and $m^+$ denote ${\rm min}(M)$ and ${\rm max}(M)$, respectively. We distinguish the following two cases.
		
		First suppose $\llbracket m^-, m^+ \rrbracket = V$, that is, $m^- =0$ and $m^+=n-1$. In this instance, for every $x \in \overline{M}$, we have $T(x,0) = T(x,n-1)$ because $M$ is a module of $T$, and hence $x \in X$ and $i_P(x) \in \{0, n-1\}$. Therefore $\overline{M} \subseteq X$ and, by the definition of a pairing, $|\overline{M}| \leq 2$ and hence $|\overline{M}|=1$ or $2$ because the module $M$ is nontrivial. Pick $x \in \overline{M}$. Since $i_P(x) = 0$ or $n-1$, then by interchanging $\underline{n}$ and $\underline{n}^{\star}$, we may assume $i_P(x) = 0$. If $x=1$, then the pairing $P$ is reducible because the block $\{0,1\}$ of $P$ would be an interval of $\underline{X}$, which is nontrivial because $|X| \geq 4$. Now suppose $x \geq 2$. If there is a vertex $v$ in $\rrbracket 0, x \llbracket \cap M$, we get $T(x,v) = 0 \neq T(x,0) =1$, contradicting that $M$ is a module of $T$. Therefore $\llbracket 1, x \rrbracket \subseteq  \overline{M}$. Since $x \geq 2$ and $|\overline{M}| \leq 2$, it follows that $x=2$ and $\overline{M} = \{1,2\}$. So $T(1, 0) = T(1,n-1)$ because $M$ is a module of $T$. Therefore, $i_P(1) = 0$ or $n-1$, and since $i_P(2) = 0$, we obtain $i_P(1) = n-1$. Now we have $\{0,2\} \in P$ and $\{1,n-1\} \in P$. Since $V \setminus \{0,1,2,n-1\}= \llbracket 3,  n-2 \rrbracket$, it follows that the interval $\llbracket 3,  n-2 \rrbracket \cap X$ of $\underline{X}$ is a (possibly empty at this stage) union of blocks of $P$. Therefore, if $\llbracket 3,  n-2 \rrbracket \cap X$ is nonempty, it is a nontrivial interval of $\underline{X}$ and hence $P$ is reducible. To conclude, we prove that $\llbracket 3,  n-2 \rrbracket \cap X$ is actually nonempty. Suppose for a contradiction that $\llbracket 3,  n-2 \rrbracket \cap X = \varnothing$. In this instance, it follows from (\ref{eqco}) that $|\llbracket 3,  n-2 \rrbracket| = 1$, that is, $n=5$. Since $\{0,2\} \in P$ and $\{1,n-1\} \in P$, it follows that $P =\{\{0,2\}, \{1,4\}\}$ and hence $T = {\rm Inv}(\underline{5}, \{\{0,2\}, \{1,4\}\}) = W_5$ (see Figure~\ref{Figure-T-U-W}). Thus $T$ is indecomposable, a contradiction. 
		
		Second suppose $\llbracket m^-, m^+ \rrbracket \varsubsetneq V$. In this instance, since $m^- < m^+$ and $\{0, n-1\} \cap \overline{\llbracket m^-, m^+ \rrbracket} \neq \varnothing$, then by (\ref{eqco}), the interval $J:=\llbracket m^-, m^+\rrbracket \cap X$ of $\underline{X}$ is a nonempty and proper subset of $X$. We will prove that $J$ is a union of blocks of $P$. Since $J$ is a nonempty and proper subset of $X$, this implies that the interval $J$ of $\underline{X}$ is nontrivial, and hence the pairing $P$ is reducible, which completes the proof. So Suppose for a contradiction that $J$ is not a union of blocks of $P$. There exist $x \in J$ and $y \in \overline{J}$ such that $\{x, y\} \in P$. Since $y \in X$ because $\{x, y\} \in P$ and since $y \in \overline{J}$, we have $y \in \overline{\llbracket m^-, m^+ \rrbracket}$.
		If $x \in \{m^-, m^+\}$, then $T(y, m^-)\neq T(y, m^+)$, which contradicts that $M$ is a module of $T$. Therefore $x \in \llbracket m^-+1, m^+-1\rrbracket$. In this instance, $T(x, m^-) \neq T(x, m^+)$. Since $M$ is a module of $T$ and $\{m^-, m^+\} \subseteq M$, it follows that $x \in M$. But $T(y, x) \neq T(y, m^-)$. Since $y \in \overline{M}$ and $\{x, m^-\} \subseteq M$, this again contradicts that $M$ is a module of $T$.
	\end{proof}

\begin{proof}[Proof of Theorem~\ref{proposition analogue}] 
	We may assume $V= \llbracket 0, n-1 \rrbracket$ and $\underline{V} = \underline{n}$, where $n \geq 5$. Recall that $Q$ is a quasi-pairing of $X:= \cup Q$. To begin, suppose that the first assertion does not hold. We have to prove that the tournaments $T$, $T- v_Q^-$, and $T- v_Q^+$ are all decomposable. 
	
	First suppose that $X$ is not a transversal of ${\rm mc}(\underline{n})$. By Fact~\ref{inter}, the tournament $T$ is decomposable. To show that $T- v_Q^-$ and $T- v_Q^+$ are also decomposable, consider $C \in {\rm mc}(\underline{n})$  such that $C \cap X = \varnothing$. Let $v = v_Q^-$ or $v_Q^+$. Since $v \in X$, we have $v \notin C$. Therefore, $C$ is a co-module of $\underline{n}-v$ (see \eqref{eq mc(n)}). Since $T-v = {\rm Inv}(\underline{n}-v, Q \setminus \{\{\hat{v}_Q, v\}\})$ and $\cup (Q \setminus \{\{ \hat{v}_Q, v\}\})$ is not a transversal of ${\rm mc}(\underline{n}-v)$ because $C \cap (\cup (Q \setminus \{\{\hat{v}_Q, v\}\})= C \cap (X \setminus \{v\}) = \varnothing$, it follows from Fact~\ref{inter} that $T-v$ is decomposable. 
	
	Second suppose that the quasi-pairing $Q$ of $\underline{X}$ is reducible. In this instance, there exists $R \varsubsetneq Q_{\rm part}$ such that $\cup R$ is a nontrivial interval of $\underline{X}$. Set $M:=\llbracket \min(\cup R), \max(\cup R) \rrbracket$. If there is $x$ in $M \cap (\cup (Q \setminus R))$, then $x \in X \setminus (\cup R)$ and $\min(\cup R) < x < \max(\cup R)$, which contradicts that $\cup R$ is an interval of $\underline{X}$. Therefore $\cup (Q \setminus R) \subseteq \overline{M}$. Since $R \subseteq \binom{M}{2}$, it follows that $Q \subseteq \binom{\overline{M}}{2} \cup \binom{M}{2}$ and hence $M$ is a module of $T$. Moreover, $|M| \geq 2$ because $|\cup R| \geq 2$, and $|\overline{M}| \geq 2$ because $\cup (Q \setminus R) \subseteq \overline{M}$ and $|\cup (Q \setminus R)| \geq 2$. Therefore, the module $M$ of $T$ is nontrivial, and hence the tournament $T$ is decomposable. Since $M$ is a module of $T$, then $M \setminus \{v_{Q}^-\}$ and $M \setminus \{v_{Q}^+\}$ are modules of $T-v_{Q}^-$ and $T-v_{Q}^+$, respectively.
	Moreover, since $2 \leq |M| \leq n-2$, and $|M| \geq 3$ if $B_Q \in R$, then the respective modules $M \setminus \{v_{Q}^-\}$ and $M \setminus \{v_{Q}^+\}$ of $T-v_{Q}^-$ and $T-v_{Q}^+$ are nontrivial. Therefore, the tournaments $T-v_{Q}^-$ and $T-v_{Q}^+$ are also decomposable.    
	
	Conversely, suppose that the quasi-pairing $Q$ of $\underline{X}$ is irreducible and that $X$ is a transversal of ${\rm mc}(\underline{n})$, where $n \geq 6$. Suppose that the tournament $T$ is decomposable. We have to prove that $T-v_{Q}^-$ or $T-v_{Q}^+$ is indecomposable. We have $|X| \geq 4$ because $X$ is a transversal of ${\rm mc}(\underline{n})$ and $\Delta(\underline{n}) = \left\lceil \frac{n+1}{2} \right\rceil \geq 4$ (see \ref{grand Delta}). Since $|X|$ is odd because $X = \cup Q$ and $Q$ is a quasi-pairing, it follows that $|X| \geq 5$. Now let $M$ be a nontrivial module of $T$, and let $m^-$ and $m^+$ denote ${\rm min}(M)$ and ${\rm max}(M)$, respectively. First note that since $X$ is a transversal of ${\rm mc}(\underline{n})$ and $\{\{0\}, \{n-1\}\} \subseteq {\rm mc}(\underline{n})$ (see (\ref{eq mc(n)})),
	\begin{equation} \label{eqifthen}
		\text{if } \llbracket m^-, m^+ \rrbracket \varsubsetneq \llbracket 0, n-1 \rrbracket, \text{ then }  X \cap \overline{\llbracket m^-, m^+ \rrbracket} \neq \varnothing.
	\end{equation}
	Given $\alpha \in \overline{M}\cap\llbracket m^-, m^+\rrbracket$, we have $\underline{n}(\alpha, m^-) \neq \underline{n}(\alpha, m^+)$, and $T(\alpha, m^-) = T(\alpha, m^+)$ because $M$ is a module of $T$. It follows that $\overline{M}\cap\llbracket m^-, m^+\rrbracket\subseteq X$ and
	\begin{equation}\label{l'intersection avec Q = un seul elmnt}
		|\{\{\alpha, m^-\}, \{\alpha, m^+\}\}\cap Q| = 1 \text{ for every } \alpha \in \overline{M}\cap\llbracket m^-, m^+\rrbracket.
	\end{equation}
	\begin{Claim}\label{claim analogue}
		\begin{sloppypar}
			For every $\alpha \in \overline{M}\cap\llbracket m^-, m^+\rrbracket$ and $\beta \in \iota_{Q}(\alpha)\cap M$, we have $\rrbracket \min(\alpha, \beta), \max(\alpha, \beta)\llbracket\cap M = \emptyset$.
		\end{sloppypar}
	\end{Claim}
	\begin{proof}[Proof of Claim~\ref{claim analogue}]
		Let $\alpha \in \overline{M}\cap\llbracket m^-, m^+\rrbracket$ and let $\beta \in \iota_{Q}(\alpha)\cap M$. By interchanging $\underline{n}$ and $\underline{n}^{\star}$, we may assume $\alpha < \beta$. We have to prove that $\rrbracket \alpha, \beta\llbracket\cap M = \emptyset$. Since $\{\alpha, \beta\} \in Q$ and for every $x \in \rrbracket \alpha, \beta \llbracket\cap M$, we have  $\underline{n}(\alpha, \beta) = \underline{n}(\alpha, x)$, and $T(\alpha, \beta) = T(\alpha, x)$ because $M$ is a module of $T$, we obtain $\rrbracket\alpha, \beta\rrbracket\cap M \subseteq \iota_{Q}(\alpha)$. Suppose for a contradiction that $\rrbracket \alpha, \beta\llbracket\cap M \neq \emptyset$. Since $\rrbracket\alpha, \beta\rrbracket\cap M \subseteq \iota_{Q}(\alpha)$ and $|\iota_{Q}(\alpha)|\leq 2$, it follows that $|\rrbracket\alpha, \beta\llbracket\cap M| = 1$ and hence $\iota_{Q}(\alpha) =  \{x, \beta\}$, where $x$ is the (unique) element of $\rrbracket\alpha, \beta\llbracket\cap M$. It follows from (\ref{l'intersection avec Q = un seul elmnt}) that $m^+ \in \{x,\beta\}$. Since $x < \beta \in M$ and $m^+ = \max(M)$, we obtain $\beta = m^+$. Thus $\hat{v}_Q = \alpha$, $v^-_Q = x$, $v^+_Q = m^+$, and $\rrbracket \hat{v}_Q, v^+_Q \llbracket \cap M = \{v^-_Q\}$. 
		If $B_Q$ is an interval of $\underline{X}$, then since $|B_Q| = 3$ and $|X| \geq 5$, the interval $B_Q$ of $\underline{X}$ is nontrivial, which contradicts that the quasi-pairing $Q$ of $\underline{X}$ is irreducible. Therefore, $B_Q$ is not an interval of $\underline{X}$. Since $\rrbracket \hat{v}_Q, v^+_Q \llbracket \cap M = \{v^-_Q\}$, it follows that $\rrbracket \hat{v}_Q, v^+_Q \llbracket \cap \overline{M} \neq \emptyset$. Therefore, since $m^+ = v^+_Q$, then by (\ref{l'intersection avec Q = un seul elmnt}) and by the definition of a quasi-pairing, $\llbracket m^-, m^+ \rrbracket \cap \overline{M} = \{\hat{v}_Q,y\}$ for some $y \in \rrbracket \hat{v}_Q, m^+ \llbracket$. Moreover $\{y, m^-\} \in Q$. If there exists $z$ in $\rrbracket m^-, \hat{v}_Q \llbracket \cap M$, we get $T(y,z) = 0 \neq T(y, m^+) =1$, which contradicts that $M$ is a module of $T$. Therefore $\rrbracket m^-, \hat{v}_Q \rrbracket \cap M = \varnothing$, and since $\llbracket \hat{v}_Q, m^+ \llbracket \cap M = \{v^-_Q\}$, we obtain $\rrbracket m^-, m^+ \llbracket \cap M = \{v^-_Q\}$. Thus $\llbracket m^-, m^+ \rrbracket \cap M = \{m^-, v^-_Q, v^+_Q\}$ and $\llbracket m^-, m^+ \rrbracket \cap \overline{M} = \{\hat{v}_Q,y\}$, and hence $\llbracket m^-, m^+ \rrbracket = \{m^-, y, \hat{v}_Q, v^-_Q, v^+_Q\}$. Thus, $\llbracket m^-, m^+ \rrbracket$ is the union of the blocks $B_Q$ and $\{y,m^-\}$ of $Q_{{\rm part}}$, and hence $\llbracket m^-, m^+ \rrbracket \subseteq X$. It follows that $\llbracket m^-, m^+ \rrbracket$ is an interval of $\underline{X}$. Moreover, since $|\llbracket m^-, m^+ \rrbracket| =5$ and $n \geq 6$, it follows from (\ref{eqifthen}) that $\llbracket m^-, m^+ \rrbracket \varsubsetneq X$. Therefore, the interval $\llbracket m^-, m^+ \rrbracket$ of $\underline{X}$ is nontrivial. Since $\llbracket m^-, m^+ \rrbracket$ is a union of blocks of $Q_{\rm part}$, this contradicts that the quasi-pairing $Q$ of $\underline{X}$ is irreducible, and completes the proof of the claim.	
	\end{proof}
	Now suppose for a contradiction that $|\overline{M}\cap\llbracket m^-, m^+\rrbracket| \geq 3$. Since $|\iota_{Q}(m^-)| + |\iota_{Q}(m^+)| \geq |\overline{M}\cap\llbracket m^-, m^+\rrbracket|$ by (\ref{l'intersection avec Q = un seul elmnt}), and $|\iota_{Q}(m^-)| + |\iota_{Q}(m^+)| \leq 3$ by the definition of a quasi-pairing, we obtain $|\overline{M}\cap\llbracket m^-, m^+\rrbracket| = 3$. Again by (\ref{l'intersection avec Q = un seul elmnt}), $\hat{v}_Q \in \{m^-,m^+\}$ and $\{v^-_Q, v^+_Q\} \subseteq \overline{M}\cap\rrbracket m^-, m^+\llbracket$. By interchanging $\underline{n}$ and $\underline{n}^{\star}$, we may assume $\hat{v}_Q = m^-$. Now let $\alpha$ be the (unique) element of $(\llbracket m^-, m^+ \rrbracket \cap \overline{M}) \setminus \{v^-_Q, v^+_Q\}$. Since $\hat{v}_Q = m^-$ and $\alpha \in (\llbracket m^-, m^+ \rrbracket \cap \overline{M}) \setminus \{v^-_Q, v^+_Q\}$, it follows from (\ref{l'intersection avec Q = un seul elmnt}) that $\{\alpha, m^+\} \in Q$. If $\alpha > v^+_Q$, then $\llbracket m^-, v_Q^+ \llbracket \cap \overline{M} = \{v_Q^-\}$, and since $\rrbracket m^-, v^+_Q \rrbracket \cap M = \emptyset$ by Claim~\ref{claim analogue}, we obtain $B_Q = \llbracket m^-, v_Q^+ \rrbracket$ and hence $B_Q$ is a nontrivial interval of $\underline{X}$, which contradicts that the quasi-pairing $Q$ of $\underline{X}$ is irreducible. Therefore $\alpha \in \rrbracket m^-, v^+_Q \llbracket \cap \overline{M}$. Since $\rrbracket m^-, v^+_Q \rrbracket \cap M = \llbracket \alpha, m^+ \llbracket \cap M = \varnothing$ by Claim~\ref{claim analogue}, it follows that $\rrbracket m^-, m^+ \llbracket \cap M = \varnothing$. Thus $\llbracket m^-, m^+ \rrbracket$ equals $\{\hat{v}_Q, \alpha, v^-_Q, v^+_Q, m^+\}$, which is the union of the blocks $B_Q$ and $\{\alpha, m^+\}$ of $Q_{{\rm part}}$, and hence $\llbracket m^-, m^+ \rrbracket$ is an interval of $\underline{X}$. Moreover, by (\ref{eqifthen}), the interval $\llbracket m^-, m^+ \rrbracket$ of $\underline{X}$ is nontrivial. This contradicts that the quasi-pairing $Q$ of $\underline{X}$ is irreducible. Therefore,
	\begin{equation} \label{inf2}
		|\overline{M}\cap\llbracket m^-, m^+\rrbracket| \leq 2.
	\end{equation}
	
	We now distinguish the following two cases. First suppose $m^-=0$ and $m^+=n-1$. Since $M \varsubsetneq \llbracket 0, n-1 \rrbracket$, $|\overline{M}| = 1$ or $2$ by (\ref{inf2}). Suppose for a contradiction that $|\overline{M}| = 2$. Set $\overline{M} := \{\alpha, \beta\}$ with $\alpha < \beta$. If $\iota_{Q}(0) = \{\alpha, \beta\}$ (resp. $\iota_{Q}(n-1) = \{\alpha, \beta\}$), then by Claim~\ref{claim analogue}, $\llbracket 0, \beta \rrbracket = B_Q$ (resp. $\llbracket \alpha, n-1 \rrbracket = B_Q$). Since $|X| \geq 5$, this implies that $B_Q$ is a nontrivial interval of $\underline{X}$, which contradicts that the quasi-pairing $Q$ of $\underline{X}$ is irreducible. Therefore $\iota_{Q}(0) \neq \{\alpha, \beta\}$ and $\iota_{Q}(n-1) \neq \{\alpha, \beta\}$. Then by (\ref{l'intersection avec Q = un seul elmnt}), $\{\{0, \alpha\}, \{n-1,\beta\}\} \subseteq Q$ or $\{\{0, \beta\}, \{n-1, \alpha\}\} \subseteq Q$. If $\{\{0, \beta\}, \{n-1, \alpha\}\} \subseteq Q$, then by Claim~ \ref{claim analogue}, $\llbracket 0, n-1 \rrbracket = \{0,n-1, \alpha, \beta\}$, contradicting $n \geq 6$. So $\{\{0, \alpha\}, \{n-1,\beta\}\} \subseteq Q$. It follows from Claim~ \ref{claim analogue} that $\alpha=1$ and $\beta = n-2$. Since $\{\{0,1\}, \{n-2,n-1\}\} \subseteq Q$, then by the definition of a quasi-pairing, $\{0,1\}$ or $\{n-2,n-1\}$ is a block of $Q_{{\rm part}}$, and hence a nontrivial interval of $\underline{X}$. This again contradicts that the quasi-pairing $Q$ of $\underline{X}$ is irreducible. Therefore $|\overline{M}| = 1$. Denote by $\alpha$ the (unique) element of  $\overline{M}$. By (\ref{l'intersection avec Q = un seul elmnt}), we have $\{0, \alpha\}\in Q$ or $\{n-1, \alpha\}\in Q$. By interchanging $\underline{n}$ and $\underline{n}^{\star}$, we may assume $\{0, \alpha\} \in Q$. It follows from  Claim~\ref{claim analogue} that $\alpha = 1$, and hence $\overline{M}=\{1\}$. Since $\{0,1\} \in Q$ and $M=\llbracket 0, n-1 \rrbracket \setminus \{1\}$ is a module of $T$, we obtain $\iota_{Q}(1) = \{0\}$. Since $Q$ is irreducible, it follows that $\iota_{Q}(0) = \{1, \beta\}$ for some $\beta \in \llbracket 3, n-1 \rrbracket$. Thus $\hat{v}_Q =0$, $v^-_Q = 1$, and $v^+_Q = \beta$. Suppose to the contrary that the tournament $T- v^-_Q$, which is $T-1$, is decomposable. We have $T-1 = \text{Inv}(\underline{n} -1, P)$, where $P := Q \setminus \{\{0,1\}\}$. Since $X$ is a transversal of ${\rm mc}(\underline{n})$ and ${\rm mc}(\underline{n}-1) = {\rm mc}(\underline{n}) \setminus \{\{1,2\}\} \subseteq {\rm mc}(\underline{n})$, then $X \setminus \{1\}$, which is $\cup P$, is also a transversal of ${\rm mc}(\underline{n}-1)$. It follows from Theorem~\ref{V souligné} that the pairing $P$ of $\underline{X}-1$ is reducible. So there exists $R \varsubsetneq P$ such that $\cup R$ is a nontrivial interval of $\underline{X} -1$. If $\{0, v^+_Q\} \notin R$, then $\cup R \subseteq \llbracket 2, n-1 \rrbracket$ and hence $\cup R$ is also a nontrivial interval of $\underline{X}$, which contradicts that $Q$ is irreducible because in this instance $R \varsubsetneq Q_{\rm part}$. If $\{0, v^+_Q\} \in R$, then $(X \setminus \{1\}) \cap \llbracket 0, v^+_Q \rrbracket  \subseteq \cup R$ and consequently $(\cup R) \cup \{1\}$, which is  $\cup ((R \setminus \{\{0,v^+_Q\}\})\cup \{B_Q\})$ and hence a union of blocks of $Q_{\rm part}$, is a nontrivial interval of $\underline{X}$, which again contradicts the irreducibility of $Q$. Therefore $T-v^-_Q$ is indecomposable. 
	
	Second suppose $\llbracket m^-, m^+\rrbracket \subsetneq \llbracket 0, n-1 \rrbracket$. We distinguish two cases. To begin, suppose that there exist $x \in M$ and $y \in \overline{\llbracket m^-, m^+ \rrbracket}$ such that $\{x, y\} \in Q$. In this instance, $\{\{y,z\}: z \in M\} \subseteq Q$ because $M$ is a module of $T$. It follows from the definition of a quasi-pairing that $M = \{m^-, m^+\}$. Moreover, $m^+ = m^-+1$ (see (\ref{l'intersection avec Q = un seul elmnt})). Thus $\hat{v}_Q=y$, $v^-_Q=m^-$, and $v^+_Q=m^-+1$. We now prove that the tournament $T-v_Q^-$ is indecomposable. Since $T-v^-_Q = \text{Inv}(\underline{n} -v^-_Q, P)$ where $P := Q \setminus \{\{\hat{v}_Q, v_Q^-\}\}$, and since $P$ is a pairing of $X \setminus \{v_Q^-\}$, then by Theorem~\ref{V souligné}, it suffices to prove that (1) $X \setminus \{v_Q^-\}$ is a transversal of ${\rm mc}(\underline{n} -v^-_Q)$, and (2) the pairing $P$ of $\underline{X} - v_Q^-$ is irreducible. For (1), let $C \in {\rm mc}(\underline{n} -v^-_Q)$. If $C \in {\rm mc}(\underline{n})$, then $X \cap C \neq \varnothing$ because $X$ is a transversal of ${\rm mc}(\underline{n})$, and since $v_Q^- \notin C$ because $C \in {\rm mc}(\underline{n} -v^-_Q)$, we obtain $(X \setminus \{v_Q^-\}) \cap C \neq \varnothing$. Now suppose $C \in {\rm mc}(\underline{n} - v_Q^-) \setminus {\rm mc}(\underline{n})$. In this instance, either $v_Q^-=0$ and $C = \{v_Q^-+1\} = \{v_Q^+\}$, or $v_Q^- \in \llbracket 2, n-3 \rrbracket$ and $C = \{v_Q^- - 1, v_Q^- +1\} = \{v_Q^- -1, v_Q^+\}$. In both instances, we have $v_Q^+ \in C$, and since $v_Q^+ \in X \setminus \{v_Q^-\}$, we obtain $(X \setminus \{v_Q^-\}) \cap C \neq \varnothing$. Therefore (1) holds. Now suppose for a contradiction that (2) does not hold. There exists $R \varsubsetneq P$ such that $\cup R$ is a nontrivial interval of $\underline{X} -v^-_Q$. If $\{\hat{v}_Q, v^+_Q\} \notin R$, then $\cup R \subseteq \llbracket 0, v^-_Q \llbracket$ or $\cup R \subseteq \rrbracket v^+_Q, n-1 \rrbracket$, and hence $\cup R$ is also a nontrivial interval of $\underline{X}$, which contradicts that $Q$ is irreducible because in this instance $R \varsubsetneq Q_{\rm part}$. If $\{\hat{v}_Q, v^+_Q\} \in R$, then $(X \setminus\{v^-_Q\}) \cap \llbracket \min(B_Q), \max(B_Q)\rrbracket \subseteq \cup R$ and consequently $(\cup R) \cup \{v^-_Q\}$, which is  $\cup ((R \setminus \{\{\hat{v}_Q,v^+_Q\}\})\cup \{B_Q\})$ and hence a union of blocks of $Q_{\rm part}$, is a nontrivial interval of $\underline{X}$, which again contradicts the irreducibility of $Q$. Therefore (2) holds.    
	
	To finish, suppose that 	
	\begin{equation}\label{{y, x} notin Q}
		Q \cap \{\{u,v\} : u \in M \text{ and } v \in \overline{\llbracket m^-, m^+ \rrbracket}\} = \varnothing.
	\end{equation}
We first show that
\begin{equation} \label{xy exist}
	\text{there exist } x \in \overline{M} \cap \llbracket m^-, m^+\rrbracket \text{ and } y \in \overline{\llbracket m^-, m^+ \rrbracket} \text{ such that } \{x, y\} \in Q.
	\end{equation}
Suppose not. In this instance, $\llbracket m^-, m^+ \rrbracket \cap X$ is a union of blocks of $Q_{\rm part}$. Since $\overline{\llbracket m^-, m^+ \rrbracket} \cap X \neq \varnothing$ (see \eqref{eqifthen}), and $\llbracket m^-, m^+ \rrbracket \cap X \neq \varnothing$ because $X$ is a transversal of ${\rm mc}(\underline{n})$, it follows that the interval $\llbracket m^-, m^+ \rrbracket \cap X$ of $\underline{X}$, which is a union of blocks of $Q_{\rm part}$,  is nontrivial. This contradicts that the quasi-pairing $Q$ of $\underline{X}$ is irreducible. Therefore \eqref{xy exist} holds.    
	By interchanging $\underline{n}$ and $\underline{n}^{\star}$, we may assume $\iota_{Q}(x) = \{m^+, y\}$ (see (\ref{l'intersection avec Q = un seul elmnt})). Thus 
	\begin{equation} \label{xvv}
	\hat{v}_Q=x \text{ and } \{v_Q^-, v_Q^+\} = \{m^+, y\}. 
	\end{equation}
	By Claim~\ref{claim analogue}, we have $\llbracket x, m^+ \llbracket \cap M = \varnothing$. Moreover $|\rrbracket x, m^+ \rrbracket \cap \overline{M}| \leq 1$ (see \eqref{inf2}). Consequently, 
	\begin{equation} \label{EQ}
		\rrbracket x, m^+ \llbracket = \varnothing, \text{ or } \rrbracket x, m^+ \llbracket = \{z\} \text{ and } \iota_Q(z) = \{m^-\} \text{ for some } z \in \overline{M} \ (\text{see} \  \eqref{l'intersection avec Q = un seul elmnt}).
	\end{equation}
We now prove that 
\begin{equation}\label{m-,x vide}
	\rrbracket m^-, x \llbracket = \varnothing, \text{ that is, } x = m^-+1.
\end{equation}
Suppose not. We distinguish two cases. First suppose $\rrbracket m^-, x \llbracket \cap \overline{M} \neq \varnothing$. Since $x \in \llbracket m^-, m^+ \rrbracket \cap \overline{M}$, so by \eqref{inf2}, $\rrbracket m^-, x \llbracket \cap \overline{M}$ is a singleton $\{m\}$. 
By \eqref{l'intersection avec Q = un seul elmnt} and \eqref{xvv}, $\{m^-, m\}\in Q_{\rm part}$. It follows from Claim~\ref{claim analogue} that $m = m^-+1$, and hence $\{m^-, m^-+1\}$ is a nontrivial interval of $\underline{X}$. Since $\{m^-, m^-+1\} \in Q_{\rm part}$, this contradicts that $Q$ is irreducible. Second, suppose $\rrbracket m^-, x \llbracket \cap \overline{M} = \varnothing$ and hence $\rrbracket m^-, x \llbracket \cap M \neq \varnothing$. So $m^-+1 \in M$. Suppose to the contrary that $\rrbracket x, m^+ \llbracket \neq \varnothing$. It follows from \eqref{EQ} that $\rrbracket x, m^+ \llbracket = \{z\}$ and $\iota_Q(z) = \{m^-\}$ for some $z \in \overline{M}$. Thus $T(z,m^-) = 1 \neq T(z, m^-+1)=0$, contradicting that $M$ is a module of $T$. So $\rrbracket x, m^+ \llbracket = \varnothing$. Therefore, it follows from \eqref{{y, x} notin Q} and \eqref{xvv} that $\llbracket m^-, x\llbracket \cap X$ is a union of blocks of $Q_{\rm part}$. Moreover, this union is nonempty because $X$ is a transversal of ${\rm mc}(\underline{n})$ and $\llbracket m^-, x\llbracket$ is a nontrivial interval of $\underline{n}$. Thus, the interval $\llbracket m^-, x\llbracket \cap X$ of $\underline{X}$ is nontrivial. Since this interval is a union of blocks of $Q_{\rm part}$, this contradicts the irreducibility of $Q$. Therefore \eqref{m-,x vide} holds. It follows from \eqref{{y, x} notin Q}, \eqref{xvv}, \eqref{EQ}, and \eqref{m-,x vide} that
	\begin{equation} \label{xz}
		\begin{cases}
			\rrbracket m^-, m^+ \llbracket = \{\hat{v}_Q\}, \ \hat{v}_Q \in \overline{M}, \ \iota_{Q}(\hat{v}_Q) = \{m^+, y\}, \text{ and } m^- \notin X; \\
			  \text{ or } \\
			\rrbracket m^-, m^+ \llbracket = \{\hat{v}_Q, z\}, \ \hat{v}_Q < z \in \overline{M}, \ \iota_{Q}(\hat{v}_Q) = \{m^+, y\}, \text{ and }\{z,m^-\} \in Q. 
		\end{cases}
	\end{equation}   
	We will now prove that the tournament $T-m^+$ is indecomposable, which completes the proof because $m^+ = v_Q^-$ or $v_Q^+$. Since $T-m^+ = \text{Inv}(\underline{n} -m^+, P)$ where $P := Q \setminus \{\{m^+, \hat{v}_Q\}\}$, and since $P$ is a pairing of $X \setminus \{m^+\}$, then by Theorem~\ref{V souligné}, it suffices to prove that (1) $X \setminus \{m^+\}$ is a transversal of ${\rm mc}(\underline{n} -m^+)$, and (2) the pairing $P$ of $\underline{X} - m^+$ is irreducible.
	To begin, let $C \in {\rm mc}(\underline{n} -m^+)$. If $C \in {\rm mc}(\underline{n})$, then $X \cap C \neq \varnothing$ because $X$ is a transversal of ${\rm mc}(\underline{n})$, and since $m^+ \notin C$ because $C \in {\rm mc}(\underline{n} -m^+)$, we obtain $(X \setminus \{m^+\}) \cap C \neq \varnothing$. Now suppose $C \in {\rm mc}(\underline{n} - m^+) \setminus {\rm mc}(\underline{n})$. In this instance, either $m^+=n-1$ and $C = \{m^+-1\} = \{n-2\}$, or $m^+ \in \llbracket 3, n-3 \rrbracket$ and $C = \{m^+ - 1, m^++1\}$. In both instances, we have $m^+-1 \in C$,  and since $m^+-1 \in X \setminus \{m^+\}$ by (\ref{xz}), we obtain $(X \setminus \{m^+\}) \cap C \neq \varnothing$. Therefore (1) holds. Now suppose for a contradiction that (2) does not hold. 
	So there exists $R \varsubsetneq P$ such that $\cup R$ is a nontrivial interval of $\underline{X} -m^+$. First suppose $\{\hat{v}_Q, y\} \notin R$. In this instance, we have $R \varsubsetneq Q_{\rm part}$. Moreover, since $\hat{v}_Q \notin \cup R$, it follows from (\ref{xz}) that $\cup R \subseteq \llbracket 0, \hat{v}_Q \llbracket$ or $\cup R \subseteq \rrbracket m^+, n-1 \rrbracket$, and hence $\cup R$ remains a nontrivial interval in $\underline{X}$. This contradicts that $Q$ is irreducible.
	Second suppose $\{\hat{v}_Q, y\} \in R$. In this instance, by using (\ref{xz}), it is straightforward to verify that $(X\setminus\{m^+\}) \cap \llbracket \min(B_Q), \max(B_Q) \rrbracket \subseteq \cup R$, and hence $(\cup R) \cup \{m^+\}$, which is $\cup ((R \setminus \{\{\hat{v}_Q, y\}\}) \cup \{B_Q\})$, is a nontrivial interval of $\underline{X}$. Since $(\cup R) \cup \{m^+\}$ is a union of blocks of $Q_{\rm part}$, this again contradicts that $Q$ is irreducible. Therefore (2) holds.
\end{proof}

\begin{proof}[Proof of Theorem~\ref{résultat plus fort}]
	First suppose that Conditions~${\rm (C1)}$--${\rm (C4)}$ are not all satisfied. We have to prove that $T$ is decomposable. If Condition~${\rm (C1)}$ is not satisfied, then $T$ is decomposable by Theorem~\ref{proposition analogue}. If Condition~${\rm (C2)}$ is not satisfied, that is, $v_Q^+ = v_Q^- +1$, then $\{v_Q^-, v_Q^+\}$ is a nontrivial module of $T$, and hence $T$ is decomposable.
	Now suppose that Condition~${\rm (C3)}$ is not satisfied, that is, there exists $v \in V$ such that $\{\{v, v+2\}, \{v+1,v+3\}\} \subseteq Q$ and $\hat{v}_Q \notin \{v, v+3\}$. In this instance, $\{v, v+3 \}$ is a nontrivial module of $T$, and hence $T$ is decomposable. Finally, suppose that Condition~${\rm (C4)}$ is not satisfied. In this instance, there exists $v \in V$ such that $\{v,v+1\} \in Q$. Moreover, $\hat{v}_Q \notin \{v,v+1\}$ or $\{\hat{v}_Q-1, \hat{v}_Q+1\} \nsubseteq \cup Q$. If $\hat{v}_Q \notin \{v,v+1\}$, then $\{v,v+1\}$ is a nontrivial module of $T$ and hence $T$ is decomposable. Now suppose $\hat{v}_Q \in \{v,v+1\}$. In this instance $\{\hat{v}_Q-1, \hat{v}_Q+1\} \nsubseteq \cup Q$. By interchanging $\underline{n}$ and $\underline{n}^{\star}$, we may assume $\hat{v}_Q-1 \notin \cup Q$. Since $\hat{v}_Q \in \{v,v+1\}$ and $\{v,v+1\} \in Q$, it follows that $\hat{v}_Q =v$. First suppose $\hat{v}_Q =v= 0$. In this instance $\{0,1\} \in Q$. Since $\hat{v}_Q \neq 1$, it follows that $T(1, V \setminus \{1\})=1$ and hence $V \setminus \{1\}$ is a nontrivial module of $T$. Thus $T$ is decomposable. Second suppose $\hat{v}_Q =v  \neq 0$. In this instance, $\{v - 1, v+1\}$ is a nontrivial module of $T$, so $T$ is decomposable.
	
Conversely, suppose that Conditions ${\rm (C1)}$--${\rm (C4)}$ are satisfied. We have to prove that the tournament $T$ is indecomposable. For $n=5$, $T$ is indecomposable as one of the eleven indecomposable tournaments shown in Figure~\ref{Figure-p-q-p}. Now suppose $n \geq 6$. Suppose to the contrary that $T$ is decomposable. Since $Q$ is an irreducible quasi-pairing of the transversal $\cup Q$ of ${\rm mc}(\underline{n})$ (see Condition~${\rm (C1)}$), it follows from Theorem~\ref{proposition analogue} that $T-v^-_Q$ or $T-v^+_Q$ is indecomposable. By interchanging $\underline{n}$ and $\underline{n}^{\star}$, we may assume that $T-v^-_Q$ is indecomposable.
Let $M$ be a nontrivial module of $T$. We distinguish the following two cases. 

First suppose $v^-_Q \notin M$. In this instance, $M$ is also a module of $T-v^-_Q$. Since $T-v^-_Q$ is indecomposable and $|M| \geq 2$, it follows that $M = V(T-v^-_Q) = \llbracket 0, n-1 \rrbracket \setminus \{v^-_Q\}$. Since $v^+_Q \in M$, $v^-_Q \notin M$, and $T(v^-_Q, v^+_Q) = 1$, this yields $T(v^-_Q, M)=1$ and hence $\hat{v}_Q = 0$ and $v^-_Q = 1$, which contradicts Condition~${\rm (C4)}$.

 Second suppose $v^-_Q \in M$. Since $M \setminus \{v^-_Q\}$ is a module of $T-v^-_Q$ and $T-v^-_Q$ is indecomposable, the module $M\setminus \{v^-_Q\}$ of $T-v^-_Q$ is trivial. Since $2 \leq |M| \leq n-1$ because $M$ is a nontrivial module of $T$, it follows that $|M \setminus \{v^-_Q\}| = 1$ and hence $M = \{v^-_Q, u\}$ for some $u \in \llbracket 0, n-1 \rrbracket \setminus\{v^-_Q\}$. Recall that by Condition~${\rm (C2)}$, $v_Q^+ > v_Q^-+1$. So $T(v_Q^-+1, v_Q^+) =1 \neq T(v_Q^-+1, v_Q^-)=0$. Since $\{v_Q^-,u\}$ is a module of $T$, it follows that $u \neq v_Q^+$. Suppose for a contradiction that $u = \hat{v}_Q$. In this instance, $T(v_Q^+, \hat{v}_Q) = T(v_Q^+, v_Q^-) =0$, and hence $\hat{v}_Q > v_Q^+$. Since $v_Q^+ > v_Q^-+1$, it follows that $T(v_Q^-+1, \hat{v}_Q) = 1 \neq T(v_Q^-+1, v_Q^-) = 0$, which contradicts that $\{v_Q^-, \hat{v}_Q\}$ is a module of $T$. Therefore $u \neq \hat{v}_Q$. Thus 
 \begin{equation} \label{equ}
 u \notin \{\hat{v}_Q, v^+_Q\}. 
 \end{equation}
  If $|v_Q^- -u| \geq 4$, then since $u \neq \hat{v}_Q$, there exists $x \in \rrbracket \min(u, v_Q^-), \max(u, v_Q^-) \llbracket$ such that $\{x,u\} \notin Q$ and $\{x,v_Q^-\} \notin Q$ and hence $T(x,u) \neq T(x,v_Q^-)$, which contradicts that $\{u, v_Q^-\}$ is a module of $T$. Therefore $|v_Q^- -u| \leq 3$. We then distinguish the following three cases.
\begin{itemize}
	\item Suppose $|v_Q^- -u| = 1$. In this instance, it follows from \eqref{equ} that $T(\hat{v}_Q, u) \neq T(\hat{v}_Q, v_Q^-)$, which contradicts that $\{u, v_Q^-\}$ is a module of $T$.
	 
	\item  Suppose $|v_Q^- -u| = 2$. Let $x:= \min(u, v_Q^-)+1$. If $x \neq \hat{v}_Q$, then it follows from \eqref{equ} that $T(\hat{v}_Q, u) \neq T(\hat{v}_Q, v_Q^-)$, contradicting that $\{u, v_Q^-\}$ is a module of $T$. Now suppose $x = \hat{v}_Q$. In this instance, it follows from Condition~${\rm (C4)}$ that $u \in \cup Q$. Moreover, $i_Q(u) \neq x$ because $u \neq v_Q^+$ (see \eqref{equ}). It follows that $T(i_Q(u), u) \neq T(i_Q(u), v_Q^-)$, which again contradicts that $\{u, v_Q^-\}$ is a module of $T$.
	
	\item Suppose $|v_Q^- -u| = 3$. Let $x$ and $y$ be the elements of $\rrbracket \min(u, v_Q^-), \max(u, v_Q^-) \llbracket$. Since $\{u, v_Q^-\}$ is a module of $T$, so necessarily $u \in \cup Q$ and $\{i_Q(u), \hat{v}_Q\} = \{x,y\}$. Thus, $|u- i_Q(u)| = |v_Q^- - \hat{v}_Q|=1$ or $|u- i_Q(u)| = |v_Q^- -\hat{v}_Q|=2$. The first instance is not possible due to Condition~${\rm (C4)}$. The second one is also not possible due to Condition~${\rm (C3)}$. \qedhere
\end{itemize}
\end{proof}


{}	
\end{document}